\documentclass[reqno]{amsart}
\usepackage{setspace}
\usepackage{graphicx}
\usepackage{amsmath}
\usepackage{amsfonts}
\usepackage{amssymb}
\begin{document}

\centerline{\huge \bf On Generalized Kissing Numbers of Convex Bodies}

\bigskip
\centerline{\large Yiming Li and Chuanming Zong}

{\large
\vspace{0.8cm}
\centerline{\begin{minipage}{12.8cm}
{\bf Abstract.} In 1694, Gregory and Newton proposed the problem to determine the kissing number of a rigid material ball. This problem and its higher dimensional generalization have been studied by many mathematicians, including Minkowski, van der Waerden, Hadwiger, Swinnerton-Dyer, Watson, Levenshtein, Odlyzko, Sloane and Musin.  In this paper, we introduce and study a further generalization of the kissing numbers for convex bodies and obtain some exact results, in particular for balls in dimensions three, four and eight.
\end{minipage}}


\vspace{1cm}
\noindent
{\LARGE\bf 1. Introduction}

\bigskip
\noindent
In 1694, Gregory and Newton discussed the following problem: {\it Can a rigid material ball be brought into contact with thirteen other such balls of the same size}? Gregory believed \lq\lq yes", while Newton thought \lq\lq no". The solution of this problem has a complicated history! Several authors claimed proofs that {\it the largest number of nonoverlapping unit balls which can be brought into contact with a fixed one is twelve} (see Hoppe \cite{Ho01}, G\"{u}nter \cite{Gu01}, Sch\"{u}tte and van der Waerden \cite{Sc01} and Leech \cite{Le02}). However, only Sch\"{u}tte and van der Waerden's proof is complete!

Let $K$ be an $n$-dimensional convex body and let $C$ be an $n$-dimensional centrally symmetric convex body centered at the origin of $\mathbb{E}^n$. let $\kappa(K)$ and $\kappa^*(K)$ denote the translative kissing number and the lattice kissing number of $K$, respectively. In other words, $\kappa(K)$ is the maximum number of nonoverlapping translates $K+{\bf x}$ that can touch $K$ at its boundary, and $\kappa^*(K)$ is defined similarly, with the restriction that the translates are members of a lattice packing of $K$. It is easy to see that
$$\kappa^*(K)\leq \kappa(K)$$
holds for every convex body $K$.

In 1904, Minkowski \cite{Mi01} defined
$$D(K)=\{{\bf x}-{\bf y}:\ {\bf x}\in K,\ {\bf y}\in K\}$$
and proved
$$\kappa^*(K)=\kappa^*(D(K))$$
and
$$\kappa(K)=\kappa(D(K)).$$
Clearly, $D(K)$ is always centrally symmetric and centered at the origin of $\mathbb{E}^n$. Three years later,  Minkowski \cite{Mi02} proved that
$$\kappa^*(K)\le 3^n-1$$
holds for every $n$-dimensional convex body $K$. In 1957, Hadwiger \cite{Ha01} improved Minkowski's upper bound to
$$\kappa(K)\le 3^n-1.$$

Let $B^n$ denote the $n$-dimensional unit ball centered at the origin of $\mathbb{E}^n$. The kissing numbers $\kappa^*(B^n)$ and $\kappa(B^n)$ have been studied by many authors (see \cite{Ba01, Le01, Mu01, Mu02, Od01, Wa01}). The known exact results are summarized in the following table.

\medskip
\centerline{
\renewcommand\arraystretch{1.5}
\begin{tabular}{|c|c|c|c|c|c|c|c|c|c|c|}
\hline $n$ & $2$ & $3$ & $4$ & $5$ & $6$ & $7$ & $8$ & $9$ & $24$\\
\hline $\kappa^*(B^n)$ & $6$ & $12$ & $24$ & $40$ & $72$ & $126$ & $240$ & $272$ & $196560$\\
\hline $\kappa(B^n)$ & $6$ & $12$ & $24$ & $??$ & $??$ & $??$ & $240$ & $??$ & $196560$\\
\hline
\end{tabular}}

\bigskip
It is well-known that each centrally symmetric convex body $C$ centered at the origin defines a metric $\|\cdot \|_C$ on $\mathbb{R}^n$ by
$$\|{\bf x}, {\bf y}\|_C=\|{\bf x}-{\bf y}\|_C=\min \{ r:\ r>0,\ {\bf x}-{\bf y}\in rC\}.$$
Especially, we use $||\cdot||$ to denote the metric defined by $B^n$.

Clearly, the kissing numbers $\kappa(C)$ and $\kappa^*(C)$ only consider the closest neighbours of $C$ in translative packings and lattice packings, respectively. In fact, in many physic situations the neighbours in a larger region also have effect on $C$. For example, in some potential energy models. Therefore, it is reasonable to make following generalizations:  For $\alpha\ge 0$, we define $\kappa_{\alpha}(C)$ to be the maximum number of translates $C+{\bf x}$ which can be packed into the region $(3+\alpha)C\setminus {\rm int}(C)$ and define $\kappa_{\alpha}^*(C)$ to be the maximum number of translates $C+{\bf x}$ which can be packed into the region $(3+\alpha)C\setminus {\rm int}(C)$ where all the translative vectors simultaneously belong to a lattice, where ${\rm int}(C)$ denotes the interior of $C$.

In this paper, among other things, we will prove the following results:

\medskip\noindent
{\bf Theorem 1.} {\it In $\mathbb{E}^2$, we have}
$$\kappa^*_\alpha(B^2)=\left\{
\begin{aligned}
6,& \ \ 0\leq\alpha<2\sqrt2-2,\\
8,& \ \ 2\sqrt2-2\leq\alpha<2\sqrt3-2,\\
12,& \ \ \alpha=2\sqrt3-2,\\
\end{aligned}
\right.
$$

\medskip\noindent
{\bf Theorem 2.} {\it In $\mathbb{E}^3$, we have}
$$\kappa^*_\alpha(B^3)=\left\{
\begin{array}{ll}
12,&  0\leq\alpha<\frac{4}{3}\sqrt3-2,\\
14,& \frac{4}{3}\sqrt3-2\leq\alpha<2\sqrt2-2,\\
20,& \alpha=2\sqrt2-2,
\end{array}
\right.
$$

\medskip\noindent
{\bf Theorem 3.} {\it In $\mathbb{E}^4$, we have}
$$\kappa^*_\alpha(B^4)=\left\{
\begin{array}{ll}
30,&  \sqrt6-2\leq\alpha<2\sqrt2-2,\\
50,& \alpha=2\sqrt2-2,
\end{array}
\right.
$$

\medskip\noindent
{\bf Theorem 4.} {\it In $\mathbb{E}^8$, when $\alpha=2\sqrt2-2$ we have}
$$\kappa^*_\alpha(B^8)=2400.$$

\vspace{1cm}
\noindent
{\LARGE\bf 2. Some Basic Lemmas}

\bigskip
In 1907, H. Minkowski \cite{Mi02} studied the lattice kissing number of an $n$-dimensional convex body and proved the following result.

\medskip\noindent
{\bf Lemma 2.1.} {\it If $K$ is an $n$-dimensional convex body, then
$$\kappa^*(K)\leq 3^n-1,$$
where the equality holds if and only if $K$ is a parallelepiped. If $C$ is an $n$-dimensional centrally symmetric strictly convex body centered at ${\bf o}$, then }
$$\kappa^*(C)\leq 2(2^n-1).$$

\medskip
For an non-negative number $\alpha$ and a packing lattice $\Lambda $ of $B^n$, we define
$$X(\alpha, \Lambda)=\{{\bf v}: 2\leq ||{\bf v}||\leq2+\alpha, {\bf v}\in\Lambda\}.$$
Next we introduce two technical results which will be frequently used in this paper.

\medskip
\noindent
{\bf Lemma 2.2.} {\it When $0\leq\alpha<2\sqrt2-2$, we have}
$$\kappa_{\alpha}^*(B^n)\leq 2(2^n-1).$$

\begin{proof} On the contrary, suppose that there are such a positive $\alpha$ which less than $2\sqrt2-2$ and a suitable lattice $\Lambda$ satisfying
$${\rm card}\{X(\alpha, \Lambda)\}\geq 2^{n+1}.$$

For convenience, we assume that ${\bf a}_1, {\bf a}_2,\ldots, {\bf a}_n$ is a basis of $\Lambda$ and say two lattice vectors
$${\bf v}=z_1{\bf a}_1+z_2{\bf a}_2+\ldots+z_n{\bf a}_n$$
and
$${\bf v}'=z_{1}'{\bf a}_1+z_2'{\bf a}_2+\ldots+z_n'{\bf a}_n$$
are equivalent if
$$z_i-z_i'\equiv 0 \quad ({\rm mod}\ 2)$$
for all $i=1,2,\ldots,n.$ In other words, ${\bf v}$ and ${\bf v}'$ are equivalent if and only if
$$\mbox{$\frac{1}{2}$}({\bf v}-{\bf v}')\in\Lambda.$$
Clearly, this relation divides the points of $\Lambda$ into $2^n$ classes.

Since $X(\alpha, \Lambda)$ is centrally symmetric and ${\rm card}\{X(\alpha, \Lambda)\}\geq 2^{n+1}$, it contains $2^n$ lattice points ${\bf v}_1, {\bf v}_2, \ldots,{\bf v}_{2^n}$ satisfying
$${\bf v}_i\neq\pm{\bf v}_j,\quad i\neq j.$$
If one of the $2^n$ points, say ${\bf v}_1$, is equivalent to ${\bf o}$, then we get
$$\mbox{$\frac{1}{2}$}{\bf v}_1\in\Lambda$$
and
$$1\leq||\mbox{$\frac{1}{2}$}{\bf v}_1||<\sqrt2<2,$$
which contradicts the assumption that $B^n+\Lambda$ is a packing.

Since all ${\bf o}$, ${\bf v}_1$, ${\bf v}_2$, $\ldots$, ${\bf v}_{2^n}$ belong to at most $2^n$ classes, two of them must belong to the same class. Without loss of generality, we may assume that ${\bf v}_1$ and ${\bf v}_2$ are equivalent. Then, we have
$$2\le ||{\bf v}_i|| < 2\sqrt2, \quad i=1,2\eqno(2.1)$$
and
$$||\mbox{$\frac{1}{2}$}({\bf v}_1\pm {\bf v}_2)||\geq2.\eqno(2.2)$$
By (2.1), we get
$$||{\bf v}_1||^{2}+||{\bf v}_2||^{2}<16.\eqno(2.3)$$
By (2.2), we get
$$||{\bf v}_1||^{2}+||{\bf v}_2||^{2}-2\langle{\bf v}_1,{\bf v}_2\rangle\geq16 \eqno(2.4)$$
and
$$||{\bf v}_1||^{2}+||{\bf v}_2||^{2}+2\langle{\bf v}_1,{\bf v}_2\rangle\geq16, \eqno(2.5)$$
where $\langle{\bf v}_1,{\bf v}_2\rangle$ denotes the inner product of ${\bf v}_1$ and ${\bf v}_2$.
Then, by (2.4) and (2.5) we obtain
$$||{\bf v}_1||^{2}+||{\bf v}_2||^{2}\geq16,$$
which contradicts (2.3). Therefore, for $0\leq\alpha<2\sqrt2-2$, we have
$$\kappa_{\alpha}^{*}(B^n)\leq 2(2^n-1).$$
Lemma 2.2 is proved.\end{proof}

\medskip
\noindent
{\bf Remark 2.1.} Writing
$$X=\{{\bf v}: 2\leq||{\bf v}|| \le 2\sqrt2,\ {\bf v}\in\Lambda\}$$
and repeating the calculations (2.1)-(2.5) one can deduce that two lattice points ${\bf v}_1, {\bf v}_2\in X$ satisfying ${\bf v}_1\neq\pm{\bf v}_2$ belong to the same equivalent class if and only if
$$||{\bf v}_1||=2\sqrt2,\quad ||{\bf v}_2||=2\sqrt2 \quad {\rm and}\quad \langle{\bf v}_1, {\bf v}_2\rangle=0.$$

\medskip
\noindent
{\bf Lemma 2.3.} {\it When $\alpha<2\sqrt3-2$,  the set $X(\alpha, \Lambda)$ contains no four collinear points.}

\begin{proof}
On the contrary, suppose $X(\alpha, \Lambda)$ has four collinear points ${\bf v}_1, {\bf v}_2, {\bf v}_3$ and ${\bf v}_4$. Without loss of generality, we may assume that $n=2$ and all ${\bf v}_1, {\bf v}_2, {\bf v}_3$ and ${\bf v}_4$ have the same $x$-coordinates, namely
${\bf v}_1=(x_0, y_1),$ ${\bf v}_2=(x_0, y_2),$ ${\bf v}_3=(x_0, y_3)$ and ${\bf v}_4=(x_0, y_4)$. Furthermore, we may also assume that
$x_0\geq 0$ and
$$y_1-y_2=y_2-y_3=y_3-y_4\ge 2,\eqno(2.6)$$
since $B^2+\Lambda$ is a packing.

If $x_0\geq\sqrt3$, we get
$$y_1-y_4\le 2\sqrt{(\alpha +2)^2-x_0^2}<2\sqrt{(2\sqrt3 )^2-3} =6 \eqno(2.7)$$
and therefore
$$\ y_1-y_2=y_2-y_3=y_3-y_4<2,$$
which contradicts (2.6).

If $x_0<\sqrt3$, we get
$$y_2-y_3\geq 2\sqrt{4-x_0^2}\eqno(2.8)$$
since both ${\bf v}_2$ and ${\bf v}_3$ belong to $\Lambda$. On the other hand, since both ${\bf v}_1$ and ${\bf v}_4$ belong to ${\rm int}(2\sqrt3 B^2)$, we get
$$y_1-y_4<2\sqrt{12-x_0^2}.\eqno(2.9)$$
By (2.8) and (2.9) one can easily deduce
$$y_1-y_4 <2\sqrt{12-x_0^2}<6\sqrt{4-x_0^2}\le 3(y_2-y_3),\eqno(2.10)$$
which contradicts (2.6).

As a conclusion of the two cases, Lemma 2.3 is proved.
\end{proof}

\medskip
\noindent
{\bf Remark 2.2.} Writing
$$X=\{{\bf v}: 2\leq||{\bf v}|| \le 2\sqrt3,\ {\bf v}\in\Lambda\}$$
and repeating the calculations (2.6)-(2.10) one can deduce that $X$ has four collinear points ${\bf v}_1, {\bf v}_2, {\bf v}_3$ and ${\bf v}_4$ if and only if, up to rotation, reflection and re-enumeration,
$${\bf v}_1=(\sqrt3, 3),\quad {\bf v}_2=(\sqrt3, 1),\quad {\bf v}_3=(\sqrt3, -1),\quad {\rm and} \quad {\bf v}_4=(\sqrt3, -3).$$

\vspace{1cm}
\noindent
{\LARGE\bf 3. Proof of Theorem 1}

\bigskip
\noindent
{\bf Theorem 1.} {\it In $\mathbb{E}^2$, we have}
$$\kappa^*_\alpha(B^2)=\left\{
\begin{aligned}
6,& \ \ 0\leq\alpha<2\sqrt2-2,\\
8,& \ \ 2\sqrt2-2\leq\alpha<2\sqrt3-2,\\
12,& \ \ \alpha=2\sqrt3-2.\\
\end{aligned}
\right.
$$

\begin{proof}
When $0\leq\alpha<2\sqrt2-2$, by Lemma 2.2 we have
$$\kappa^*_\alpha(B^2)\leq6.$$
Combining with
$$\kappa^*_\alpha(B^2)\geq\kappa^*(B^2)=6,$$
for $0\leq\alpha<2\sqrt2-2$ we get
$$\kappa^*_\alpha(B^2)=6.\eqno(3.1)$$

By Remark 2.1, when $\alpha=2\sqrt2-2$ one can deduce that
$$\kappa^*_\alpha(B^2)=8,\eqno(3.2)$$
where the equality holds when the corresponding lattice $\Lambda$ is generated by ${\bf a}_1=(2,0)$, ${\bf a}_2=(0,2)$. In fact, the optimal lattice is unique up to rotation and reflection.

When $2\sqrt2-2\leq\alpha<2\sqrt3-2$, we assume that $B^2+\Lambda$ is a lattice packing attaining $\kappa^*_\alpha (B^2)$. Then we have
$${\rm card}\{X(\alpha, \Lambda)\}\geq8.$$
Without loss of generality, by a routine argument we may assume that $\{{\bf a}_1, {\bf a}_2\}$ is a basis of $\Lambda$, $\| {\bf a}_1\|=2$ and $\{{\bf a}_1, {\bf a}_2\}\subset X(\alpha, \Lambda)$.

For an arbitrary vector ${\bf v}\in X(\alpha, \Lambda)$ which is not $\pm{\bf a}_1$, the lattice $\Lambda'$ generated by $\{{\bf a}_1, {\bf v}\}$ is a sublattice of $\Lambda$. Therefore, let $\det (\Lambda)$ denote the determinant of the lattice $\Lambda$, we have
$$\det(\Lambda')=g\det(\Lambda),$$
where $g$ is a positive integer. It is easy to see that
$$\det (\Lambda')\leq ||{\bf a}_1||\cdot||{\bf v}||<2\cdot2\sqrt3=4\sqrt3 \eqno(3.3)$$
and
$$\det(\Lambda)\geq\frac{\omega (B^2)}{\delta^*(B^2)}=2\sqrt3, \eqno(3.4)$$
where $\omega(B^2)$ denotes the area of $B^2$ and $\delta^*(B^2)=\pi/\sqrt{12}$ is the density of the densest lattice packing of $B^2$. By (3.3) and (3.4) one can easily deduce that $g=1$. Consequently, if
$${\bf v}=z_1{\bf a}_1+z_2{\bf a}_2,$$
then we must have
$$z_2=\pm g=\pm1.\eqno(3.5)$$

Since $\pm{\bf a}_2\in X(\alpha, \Lambda)$, by Lemma 2.3 one can deduce that
$$|z_1|\leq2,$$
which means that
$${\bf v}=\pm{\bf a}_2,\  \pm({\bf a}_1+{\bf a}_2),\ \pm({\bf a}_1-{\bf a}_2),\ \pm(2{\bf a}_1+{\bf a}_2),\ {\rm or}\ \pm(2{\bf a}_1-{\bf a}_2).$$
If both $({\bf a}_1+{\bf a}_2)$ and $({\bf a}_1-{\bf a}_2)$ belong to $X(\alpha, \Lambda),$ since $2{\bf a}_1+{\bf a}_2,\ {\bf a}_1+{\bf a}_2,\ {\bf a}_2,\ -{\bf a}_1+{\bf a}_2$ are collinear, $(2{\bf a}_1+{\bf a}_2)$ can not belong to $X(\alpha, \Lambda).$ Similarly, since $2{\bf a}_1-{\bf a}_2,\ {\bf a}_1-{\bf a}_2,\ -{\bf a}_2,\ -{\bf a}_1-{\bf a}_2$ are collinear, $(2{\bf a}_1-{\bf a}_2)\notin X(\alpha, \Lambda).$
On the other hand, if both $(2{\bf a}_1+{\bf a}_2)$ and $(2{\bf a}_1-{\bf a}_2)$ belong to $X(\alpha, \Lambda)$, by convexity one can deduce that
$$2{\bf a}_1\in X(\alpha, \Lambda),$$
which contradicts to $X(\alpha, \Lambda)\subset {\rm int}(2\sqrt3B^2)$. As a conclusion of these two cases, we get
$${\rm card}\{X(\alpha, \Lambda)\}\leq8.$$
Combining with (3.2), for $2\sqrt2-2\leq\alpha<2\sqrt3-2$ we get
$$\kappa_\alpha^*(B^2)=8.\eqno (3.6)$$

Finally, we deal with the case $\alpha=2\sqrt3-2$. If $X(\alpha, \Lambda)$ contains no four collinear points and (3.5) holds, by previous arguments we still obtain
$${\rm card}\{X(\alpha, \Lambda)\}\leq8.$$
Therefore, the necessary condition for ${\rm card}\{X(\alpha, \Lambda)\}\geq10$ is either $X(\alpha, \Lambda)$ contains four collinear points or (3.5) does not hold.

If $X(\alpha, \Lambda)$ has four collinear vectors ${\bf v}_1, {\bf v}_2, {\bf v}_3$ and ${\bf v}_4$, by Remark 2.2 we may assume that
$${\bf v}_1=(\sqrt3, 3),\quad {\bf v}_2=(\sqrt3, 1),\quad {\bf v}_3=(\sqrt3, -1),\quad {\bf v}_4=(\sqrt3, -3).$$
In this case, it is easy to verify that $\Lambda$ is generated by ${\bf a}_1=(\sqrt3, 1), {\bf a}_2=(\sqrt3, -1)$ and
$${\rm card}\{X(\alpha, \Lambda)\}=12.$$

If there is a point
$${\bf v}=z_1{\bf a}_1+z_2{\bf a}_2\in X(\alpha, \Lambda)$$
with $z_2=\pm2,$ by repeating (3.3) and (3.4) we get
$$\det(\Lambda)=\frac{\omega(B^2)}{\delta^*(B^2)}=2\sqrt3.$$
In this case, $\Lambda$ is the densest packing lattice of $B^2$ and therefore
$${\rm card}\{X(\alpha, \Lambda)\}=12.$$

As a conclusion of the two cases, for $\alpha=2\sqrt3-2$ we have
$$\kappa_\alpha^*(B^2)=12,\eqno(3.7)$$
and the corresponding lattice $\Lambda$ is generated by ${\bf a}_1=(\sqrt3, 1)$ and ${\bf a}_2=(\sqrt3, -1)$, up to rotation and reflection.

The theorem follows from (3.1), (3.6) and (3.7).
\end{proof}

\bigskip
In 2003, C. Zong \cite {Zo07} proved the following result.

\medskip
\noindent
{\bf Lemma 3.1.} {\it For every two-dimensional centrally symmetric convex domain $C$ there is a parallelogram with vertices ${\bf o}$, ${\bf v}_1$, ${\bf v}_2$ and ${\bf v}_1+{\bf v}_2$ such that}
$$||{\bf o}, {\bf v}_1||_C=||{\bf o}, {\bf v}_2||_C=2$$
{\it and}
$$2\leq||{\bf v}_1, {\bf v}_2||_C=||{\bf o}, {\bf v}_1+{\bf v}_2||_C\leq2\sqrt2.$$

\medskip
This result has the following corollary.

\medskip
\noindent
{\bf Corollary 3.1.} {\it When $\alpha =2\sqrt2-2,$
$$\kappa^*_\alpha (C)\ge 8$$
holds for every two-dimensional centrally symmetric convex domain $C$.}

\vspace{1cm}
\noindent
{\LARGE\bf 4. Proof of Theorem 2}

\bigskip
\noindent
{\bf Lemma 4.1.} {\it When $\frac{4}{3}\sqrt3-2\leq\alpha<2\sqrt2-2$, we have $\kappa_{\alpha}^{*}(B^3)=14$.}

\begin{proof}
Let $\Lambda$ be the lattice generated by ${\bf a}_1=(-\frac{2}{3}\sqrt6, \frac{2}{3}\sqrt3,0), {\bf a}_2=(\frac{2}{3}\sqrt6, \frac{2}{3}\sqrt3,0)$ and ${\bf a}_3=(0,\frac{2}{3}\sqrt3, \frac{2}{3}\sqrt6)$. When $\alpha=\frac{4}{3}\sqrt3-2,$ one can verify that
$${\rm card}\{X(\alpha,\Lambda)\}=14.$$
Combining with Lemma 2.2, for $\frac{4}{3}\sqrt3-2\leq\alpha<2\sqrt2-2$ we have
$$\kappa_{\alpha}^{*}(B^3)=14.$$
Lemma 4.1 is proved.
\end{proof}

\medskip
\noindent
{\bf Lemma 4.2.} {\it Let ${\bf v}_1, {\bf v}_2, {\bf v}_3$ be three linearly independent vectors of $X(\alpha,\Lambda)$ and let $\Lambda'$ denote the lattice generated by them. If $0\le \alpha<\frac{4}{3}\sqrt3-2$, then}
$${\rm det}(\Lambda')/{\rm det}(\Lambda)\le 2.$$

\begin{proof}
Since $\alpha<\frac{4}{3}\sqrt3-2$, we have
$$\det(\Lambda')\leq||{\bf v}_1||\cdot||{\bf v}_2||\cdot||{\bf v}_3||<(\mbox{$\frac{4}{3}$}\sqrt3)^3.$$
On the other hand, let $\delta^*(B^3)$ denote the density of the densest three-dimensional lattice sphere packing which is $\pi/\sqrt{18}$, we have
$$\det(\Lambda)\geq\frac{{\rm vol}(B^3)}{\delta^*(B^3)}=\frac{\frac{4}{3}\pi}{\frac{\pi}{\sqrt{18}}}=4\sqrt2.$$
Therefore, we get
$${\rm det}(\Lambda')/{\rm det}(\Lambda)< (\mbox{$\frac{4}{3}$}\sqrt3)^3/4\sqrt2< 3.$$
In other words, ${\rm det}(\Lambda')/{\rm det}(\Lambda)$ only can take two values, one or two. Lemma 4.2 is proved.
\end{proof}

\medskip
\noindent
{\bf Lemma 4.3.} {\it If $\alpha < \frac{4}{3}\sqrt3-2$ and ${\rm card}\{X(\alpha, \Lambda)\}=14$, then $\Lambda$ has a basis $\{{\bf a}_1, {\bf a}_2, {\bf a}_3\}$ such that
$$X(\alpha, \Lambda)=\{\pm{\bf a}_1,\ \pm{\bf a}_2,\ \pm{\bf a}_3,\ \pm({\bf a}_1+{\bf a}_2),\ \pm({\bf a}_1+{\bf a}_3),\ \pm({\bf a}_2+{\bf a}_3),\ \pm({\bf a}_1+{\bf a}_2+{\bf a}_3)\}.$$}

\begin{proof} Suppose
$$X(\alpha,\Lambda)=\{\pm{\bf v}_1,\ \pm{\bf v}_2,\ \pm{\bf v}_3,\ \pm{\bf v}_4,\ \pm{\bf v}_5,\ \pm{\bf v}_6,\ \pm{\bf v}_7\}.$$
By Remark 2.1, ${\bf v}_i$ cannot equivalent to $\pm{\bf v}_j$ for $i\neq j$, since $\alpha < \frac{4}{3}\sqrt3-2$.

\noindent
Therefore, for a basis $\{{\bf a}_1, {\bf a}_2, {\bf a}_3\}$ of $\Lambda$, we may assume
$${\bf v}_1,\ {\bf v}_2,\ {\bf v}_3,\ {\bf v}_4,\ {\bf v}_5,\ {\bf v}_6,\ {\bf v}_7$$
were equivalent to
$${\bf a}_1,\ {\bf a}_2,\ {\bf a}_3,\ {\bf a}_1+{\bf a}_2,\ {\bf a}_1+{\bf a}_3,\ {\bf a}_2+{\bf a}_3,\ {\bf a}_1+{\bf a}_2+{\bf a}_3,$$
respectively.

It is easy to verify that, we can expand to a basis of $\Lambda$ based on \{${\bf v}_1, {\bf v}_2$\}. Without loss of generality, we assume
$${\bf v}_1=\pm{\bf a}_1,\ {\bf v}_2=\pm{\bf a}_2.$$
By Lemma 4.2, we have
$${\bf v}_3=z_1{\bf a}_1+z_2{\bf a}_2\pm{\bf a}_3,$$
where $z_1, z_2$ is even. Since \{${\bf a}_1,\ {\bf a}_2,\ z_1{\bf a}_1+z_2{\bf a}_2\pm{\bf a}_3$\} is also a basis of $\Lambda$, we assume
$${\bf v}_3=\pm{\bf a}_3$$
without loss of generality.

\noindent
Therefore, for
$${\bf v}=z_1{\bf a}_1+z_2{\bf a}_2+z_3{\bf a}_3$$
belong to $X(\alpha,\Lambda)$, we have
$$||z_i||\leq2$$
holds for $i=1,2,3$, by Lemma 4.2. Thus, we have
$${\bf v}_4=\pm{\bf a}_1\pm{\bf a}_2+z{\bf a}_3,\quad z=0\ {\rm or}\ \pm2,$$
$${\bf v}_5=\pm{\bf a}_1+z'{\bf a}_2\pm{\bf a}_3,\quad z'=0\ {\rm or}\ \pm2,$$
$${\bf v}_6=z''{\bf a}_1\pm{\bf a}_2\pm{\bf a}_3,\quad z''=0\ {\rm or}\ \pm2,$$
and
$${\bf v}_7=\pm{\bf a}_1\pm{\bf a}_2\pm{\bf a}_3.$$

Suppose one of $z, z', z''$ is $\pm2$, without loss of generality, say $z=\pm2$. Furthermore, we may assume
$${\bf v}_4=\pm({\bf a}_1+{\bf a}_2+2{\bf a}_3),$$
since the sign of ${\bf a}_1, {\bf a}_2, {\bf a}_3$ do not change $\pm{\bf v}_1, \pm{\bf v}_2, \pm{\bf v}_3$.

It is easy to verify that
$$\det({\bf a}_1-{\bf a}_3, {\bf a}_2, {\bf a}_1+{\bf a}_2+2{\bf a}_3)=\left |\begin{array}{cccc}
1 & 0 & -1\\
0 & 1 & 0\\
1 & 1 & 2\\
\end{array}\right|\det({\bf a}_1,{\bf a}_2,{\bf a}_3)=3\det(\Lambda),$$
$$\det({\bf a}_1, {\bf a}_1+2{\bf a}_2+{\bf a}_3, {\bf a}_1+{\bf a}_2+2{\bf a}_3)=\left |\begin{array}{cccc}
1 & 0 & 0\\
1 & 2 & 1\\
1 & 1 & 2\\
\end{array}\right|\det({\bf a}_1,{\bf a}_2,{\bf a}_3)=3\det(\Lambda),$$
$$\det({\bf a}_1, {\bf a}_1+2{\bf a}_2-{\bf a}_3, {\bf a}_1+{\bf a}_2+2{\bf a}_3)=\left |\begin{array}{cccc}
1 & 0 & 0\\
1 & 2 & -1\\
1 & 1 & 2\\
\end{array}\right|\det({\bf a}_1,{\bf a}_2,{\bf a}_3)=5\det(\Lambda),$$
$$\det({\bf a}_1, -{\bf a}_1+2{\bf a}_2+{\bf a}_3, {\bf a}_1+{\bf a}_2+2{\bf a}_3)=\left |\begin{array}{cccc}
1 & 0 & 0\\
-1 & 2 & 1\\
1 & 1 & 2\\
\end{array}\right|\det({\bf a}_1,{\bf a}_2,{\bf a}_3)=3\det(\Lambda),$$
$$\det({\bf a}_1, -{\bf a}_1+2{\bf a}_2-{\bf a}_3, {\bf a}_1+{\bf a}_2+2{\bf a}_3)=\left |\begin{array}{cccc}
1 & 0 & 0\\
-1 & 2 & -1\\
1 & 1 & 2\\
\end{array}\right|\det({\bf a}_1,{\bf a}_2,{\bf a}_3)=5\det(\Lambda),$$
where $\det({\bf b}_1,{\bf b}_2,{\bf b}_3)$ denotes the determinant of the lattice which generated by $\{{\bf b}_1,{\bf b}_2,{\bf b}_3\}$. By Lemma 4.2, ${\bf v}_5$ cannot be $\pm({\bf a}_1-{\bf a}_3), \pm({\bf a}_1+2{\bf a}_2+{\bf a}_3), \pm({\bf a}_1+2{\bf a}_2-{\bf a}_3), \pm(-{\bf a}_1+2{\bf a}_2+{\bf a}_3), \pm(-{\bf a}_1+2{\bf a}_2-{\bf a}_3)$. Therefore,
$${\bf v}_5=\pm({\bf a}_1+{\bf a}_3).$$
By the same deduction, we have
$${\bf v}_6=\pm({\bf a}_2+{\bf a}_3).$$

For ${\bf v}_7$, since
$$\det({\bf a}_1, {\bf a}_1+{\bf a}_2-{\bf a}_3, {\bf a}_1+{\bf a}_2+2{\bf a}_3)=\left |\begin{array}{cccc}
1 & 0 & 0\\
1 & 1 & -1\\
1 & 1 & 2\\
\end{array}\right|\det({\bf a}_1,{\bf a}_2,{\bf a}_3)=3\det(\Lambda),$$
$$\det(-{\bf a}_1, {\bf a}_1-{\bf a}_2+{\bf a}_3, {\bf a}_1+{\bf a}_2+2{\bf a}_3)=\left |\begin{array}{cccc}
-1 & 0 & 0\\
1 & -1 & 1\\
1 & 1 & 2\\
\end{array}\right|\det({\bf a}_1,{\bf a}_2,{\bf a}_3)=3\det(\Lambda),$$
$$\det(-{\bf a}_2, {\bf a}_1-{\bf a}_2-{\bf a}_3, {\bf a}_1+{\bf a}_2+2{\bf a}_3)=\left |\begin{array}{cccc}
0 & -1 & 0\\
1 & -1 & -1\\
1 & 1 & 2\\
\end{array}\right|\det({\bf a}_1,{\bf a}_2,{\bf a}_3)=3\det(\Lambda),$$
${\bf v}_7$ cannot be $\pm({\bf a}_1+{\bf a}_2-{\bf a}_3), \pm({\bf a}_1-{\bf a}_2+{\bf a}_3), \pm({\bf a}_1-{\bf a}_2-{\bf a}_3)$, by Lemma 4.2. Therefore,
$${\bf v}_7=\pm({\bf a}_1+{\bf a}_2+{\bf a}_3),$$
which means that
$$X(\alpha, \Lambda)=\{\pm{\bf a}_1, \pm{\bf a}_2, \pm{\bf a}_3, \pm({\bf a}_1+{\bf a}_2+2{\bf a}_3), \pm({\bf a}_1+{\bf a}_3), \pm({\bf a}_2+{\bf a}_3), \pm({\bf a}_1+{\bf a}_2+{\bf a}_3)\}.$$
By choose ${\bf a}_1'=-({\bf a}_1+{\bf a}_3),\ {\bf a}_2'={\bf a}_3,\ {\bf a}_3'={\bf a}_1+{\bf a}_2+{\bf a}_3$, one can verify that Lemma 4.3 holds in this case.

On the other hand, if $z=z'=z''=0$, we have
$${\bf v}_4=\pm({\bf a}_1+{\bf a}_2)\ {\rm or}\ \pm({\bf a}_1-{\bf a}_2),$$
$${\bf v}_5=\pm({\bf a}_1+{\bf a}_3)\ {\rm or}\ \pm({\bf a}_1-{\bf a}_3),$$
$${\bf v}_6=\pm({\bf a}_2+{\bf a}_3)\ {\rm or}\ \pm({\bf a}_2-{\bf a}_3).$$
Since the sign of ${\bf a}_2$ do not change $\pm{\bf v}_1, \pm{\bf v}_2, \pm{\bf v}_3$, we assume
$${\bf v}_4=\pm({\bf a}_1+{\bf a}_2),$$
without loss of generality. Furthermore, since the sign of ${\bf a}_3$ do not change $\pm{\bf v}_1, \pm{\bf v}_2, \pm{\bf v}_3,$ $\pm{\bf v}_4$, we may assume
$${\bf v}_5=\pm({\bf a}_1+{\bf a}_3),\ {\bf v}_6=\pm({\bf a}_2+{\bf a}_3)$$
or
$${\bf v}_5=\pm({\bf a}_1+{\bf a}_3),\ {\bf v}_6=\pm({\bf a}_2-{\bf a}_3).$$
Combining with
$${\bf v}_7=\pm({\bf a}_1+{\bf a}_2+{\bf a}_3),\ \pm({\bf a}_1+{\bf a}_2-{\bf a}_3), \  \pm({\bf a}_1-{\bf a}_2+{\bf a}_3) \ {\rm or}\ \pm({\bf a}_1-{\bf a}_2-{\bf a}_3),$$
we obtain that in this case, $X(\alpha, \Lambda)$ is one of the following sets:

(1). $\{\pm{\bf a}_1, \pm{\bf a}_2, \pm{\bf a}_3, \pm({\bf a}_1+{\bf a}_2), \pm({\bf a}_1+{\bf a}_3), \pm({\bf a}_2+{\bf a}_3), \pm({\bf a}_1+{\bf a}_2+{\bf a}_3)\};$

(2). $\{\pm{\bf a}_1, \pm{\bf a}_2, \pm{\bf a}_3, \pm({\bf a}_1+{\bf a}_2), \pm({\bf a}_1+{\bf a}_3), \pm({\bf a}_2+{\bf a}_3), \pm({\bf a}_1+{\bf a}_2-{\bf a}_3)\};$

(3). $\{\pm{\bf a}_1, \pm{\bf a}_2, \pm{\bf a}_3, \pm({\bf a}_1+{\bf a}_2), \pm({\bf a}_1+{\bf a}_3), \pm({\bf a}_2+{\bf a}_3), \pm({\bf a}_1-{\bf a}_2+{\bf a}_3)\};$

(4). $\{\pm{\bf a}_1, \pm{\bf a}_2, \pm{\bf a}_3, \pm({\bf a}_1+{\bf a}_2), \pm({\bf a}_1+{\bf a}_3), \pm({\bf a}_2+{\bf a}_3), \pm({\bf a}_1-{\bf a}_2-{\bf a}_3)\};$

(5). $\{\pm{\bf a}_1, \pm{\bf a}_2, \pm{\bf a}_3, \pm({\bf a}_1+{\bf a}_2), \pm({\bf a}_1+{\bf a}_3), \pm({\bf a}_2-{\bf a}_3), \pm({\bf a}_1+{\bf a}_2+{\bf a}_3)\};$

(6). $\{\pm{\bf a}_1, \pm{\bf a}_2, \pm{\bf a}_3, \pm({\bf a}_1+{\bf a}_2), \pm({\bf a}_1+{\bf a}_3), \pm({\bf a}_2-{\bf a}_3), \pm({\bf a}_1+{\bf a}_2-{\bf a}_3)\};$

(7). $\{\pm{\bf a}_1, \pm{\bf a}_2, \pm{\bf a}_3, \pm({\bf a}_1+{\bf a}_2), \pm({\bf a}_1+{\bf a}_3), \pm({\bf a}_2-{\bf a}_3), \pm({\bf a}_1-{\bf a}_2+{\bf a}_3)\};$

(8). $\{\pm{\bf a}_1, \pm{\bf a}_2, \pm{\bf a}_3, \pm({\bf a}_1+{\bf a}_2), \pm({\bf a}_1+{\bf a}_3), \pm({\bf a}_2-{\bf a}_3), \pm({\bf a}_1-{\bf a}_2-{\bf a}_3)\}.$

Since
$$\det({\bf a}_1+{\bf a}_3, {\bf a}_2+{\bf a}_3, -({\bf a}_1+{\bf a}_2-{\bf a}_3))=\left |\begin{array}{cccc}
1 & 0 & 1\\
0 & 1 & 1\\
-1 & -1 & 1\\
\end{array}\right|\det({\bf a}_1,{\bf a}_2,{\bf a}_3)=3\det(\Lambda),$$
$$\det({\bf a}_1+{\bf a}_2, {\bf a}_2+{\bf a}_3, {\bf a}_1-{\bf a}_2+{\bf a}_3)=\left |\begin{array}{cccc}
1 & 1 & 0\\
0 & 1 & 1\\
1 & -1 & 1\\
\end{array}\right|\det({\bf a}_1,{\bf a}_2,{\bf a}_3)=3\det(\Lambda),$$
$$\det({\bf a}_1+{\bf a}_2, {\bf a}_1+{\bf a}_3, {\bf a}_1-{\bf a}_2-{\bf a}_3)=\left |\begin{array}{cccc}
1 & 1 & 0\\
1 & 0 & 1\\
1 & -1 & -1\\
\end{array}\right|\det({\bf a}_1,{\bf a}_2,{\bf a}_3)=3\det(\Lambda),$$
$$\det({\bf a}_2-{\bf a}_3, {\bf a}_1+{\bf a}_3, {\bf a}_1-{\bf a}_2-{\bf a}_3)=\left |\begin{array}{cccc}
0 & 1 & -1\\
1 & 0 & 1\\
1 & -1 & -1\\
\end{array}\right|\det({\bf a}_1,{\bf a}_2,{\bf a}_3)=3\det(\Lambda),$$
by Lemma 4.2, $X(\alpha, \Lambda)$ cannot be the sets (2), (3), (4), (8).

For set (1), Lemma 4.3 already holds; For set (5), by choose ${\bf a}_1'={\bf a}_1+{\bf a}_2,\ {\bf a}_2'=-{\bf a}_2,\ {\bf a}_3'={\bf a}_3,$ it can be verified that Lemma 4.3 holds; For set (6), by choose ${\bf a}_1'={\bf a}_1,\ {\bf a}_2'={\bf a}_2-{\bf a}_3,\ {\bf a}_3'={\bf a}_3$, it can be verified that Lemma 4.3 holds; For set (7), by choose ${\bf a}_1'={\bf a}_1,\ {\bf a}_2'={\bf a}_2,\ {\bf a}_3'=-({\bf a}_1+{\bf a}_3)$, it can be verified that Lemma 4.3 holds.

As a conclusion of two cases, Lemma 4.3 is proved.
\end{proof}

\medskip
\noindent
{\bf Lemma 4.4.} {\it When $0\leq\alpha<\frac{4}{3}\sqrt3-2$, we have $\kappa_{\alpha}^{*}(B^3)=12$.}
\begin{proof}On the contrary, suppose there exist a packing lattice $\Lambda$ of $B^3$ and $0\leq\alpha<\frac{4}{3}\sqrt3-2$ satisfy
$${\rm card}\{X(\alpha,\Lambda)\}=14.$$
By Lemma 4.3, there exist a basis \{${\bf a}_1, {\bf a}_2, {\bf a}_3$\} of $\Lambda$ such that
$$X(\alpha, \Lambda)=\{\pm{\bf a}_1,\ \pm{\bf a}_2,\ \pm{\bf a}_3,\ \pm({\bf a}_1+{\bf a}_2),\ \pm({\bf a}_1+{\bf a}_3),\ \pm({\bf a}_2+{\bf a}_3),\ \pm({\bf a}_1+{\bf a}_2+{\bf a}_3)\}.$$

Without loss of generality, we suppose ${\bf a}_1,\ {\bf a}_2,\ {\bf a}_1+{\bf a}_2$ lies in plane
$$\{(v_1,v_2,v_3):\ v_3=0\}$$
and
$${\bf a}_1+{\bf a}_2=(0,c,0),\quad {\bf a}_1=(-a,b,0),\quad {\bf a}_2=(a,c-b,0).$$
Then we have:
$$\frac{16}{3}>||{\bf a}_1+{\bf a}_2||^2=c^2\geq4, \eqno(4.1)$$
$$\frac{16}{3}>||{\bf a}_1||^2=a^2+b^2\geq4, \eqno(4.2)$$
$$\frac{16}{3}>||{\bf a}_2||^2=a^2+b^2+c^2-2bc\geq4. \eqno(4.3)$$

Let ${\bf a}_3=(v_1,v_2,v_3)$ and denote $||{\bf a}_3||^2=D_1$, $||{\bf a}_1+{\bf a}_2+{\bf a}_3||^2=D_2$, then we have
$$\frac{16}{3}>D_1=v_1^2+v_2^2+v_3^2\geq4, \eqno(4.4)$$
$$\frac{16}{3}>D_2=v_1^2+v_2^2+v_3^2+c^2+2cv_2\geq4. \eqno(4.5)$$
By (4.4) and (4.5), we have
$$v_2=\frac{D_2-D_1-c^2}{2c}.\eqno(4.6)$$
By (4.4) and (4.6), we have
\begin{eqnarray*}
\frac{16}{3}>||{\bf a}_1+{\bf a}_3||^2\hspace{-0.2cm} &= \hspace{-0.2cm}&(v_1-a)^2+(v_2+b)^2+v_3^2\\
\hspace{-0.2cm} &= \hspace{-0.2cm}&D_1+a^2+b^2-2av_1+2bv_2
\end{eqnarray*}
and
\begin{eqnarray*}
\frac{16}{3}>||{\bf a}_2+{\bf a}_3||^2\hspace{-0.2cm} &= \hspace{-0.2cm}&(v_1+a)^2+(v_2+c-b)^2+v_3^2\\
\hspace{-0.2cm} &= \hspace{-0.2cm}&D_1+a^2+2av_1+2v_2(c-2b)+2bv_2+b^2+c^2-2bc\\
\hspace{-0.2cm} &= \hspace{-0.2cm}&D_1+a^2+b^2+2av_1+2bv_2+(2v_2+c)(c-2b)\\
\hspace{-0.2cm} &= \hspace{-0.2cm}&D_1+a^2+b^2+2av_1+2bv_2+\frac{D_2-D_1}{c}(c-2b).
\end{eqnarray*}
Using (4.6) again, we obtain
\begin{eqnarray*}
\frac{16}{3}\hspace{-0.2cm} &> \hspace{-0.2cm}&\frac{||{\bf a}_1+{\bf a}_3||^2+||{\bf a}_2+{\bf a}_3||^2}{2}\\
\hspace{-0.2cm} &= \hspace{-0.2cm}&D_1+a^2+b^2+2bv_2+\frac{D_2-D_1}{2c}(c-2b)\\
\hspace{-0.2cm} &= \hspace{-0.2cm}&D_1+a^2+b^2+\frac{b}{c}(D_2-D_1)-bc+\frac{D_2-D_1}{2}-\frac{b}{c}(D_2-D_1)\\
\hspace{-0.2cm} &= \hspace{-0.2cm}&\frac{D_1+D_2}{2}+a^2+b^2-bc.
\end{eqnarray*}
$$\eqno(4.7)$$

On the other hand, by (4.2) and (4.3), we have
$$2(a^2+b^2-bc)+c^2\geq8.$$
Therefore, by (4.1), we get
$$a^2+b^2-bc\geq4-\frac{c^2}{2}>\frac{4}{3},\eqno(4.8)$$
together with (4.4) and (4.5), we have
$$\frac{D_1+D_2}{2}+a^2+b^2-bc>4+\frac{4}{3}=\frac{16}{3},$$
which contradicts to (4.7).

Therefore, when $0\leq\alpha<\frac{4}{3}\sqrt3-2$, we have
$$\kappa^*_\alpha(B^3)\leq12.$$
Combining with
$$\kappa^*_\alpha(B^3)\geq\kappa^*(B^3)=12,$$
for $0\leq\alpha<\frac{4}{3}\sqrt3-2$ we have
$$\kappa^*_\alpha(B^3)=12.$$
Lemma 4.4 is proved.
\end{proof}

\medskip
\noindent
{\bf Remark 4.1.} For $\alpha=\frac{4}{3}\sqrt3-2$, by repeating (4.1)-(4.8) one can deduce that
$${\rm card}\{X(\alpha,\Lambda)\}=14$$
if and only if the lattice $\Lambda$ is generated by ${\bf a}_1=(-\frac{2}{3}\sqrt6, \frac{2}{3}\sqrt3,0), {\bf a}_2=(\frac{2}{3}\sqrt6, \frac{2}{3}\sqrt3,0)$ and ${\bf a}_3=(0,\frac{2}{3}\sqrt3, \frac{2}{3}\sqrt6)$, up to rotation and reflection.

\medskip
\noindent
{\bf Lemma 4.5.} {\it When $\alpha=2\sqrt2-2$, we have $\kappa^*_\alpha(B^3)=20$.}
\begin{proof}
Let $\Lambda$ be the lattice generated by ${\bf a}_1=(2,0,0), {\bf a}_2=(0,2,0)$ and ${\bf a}_3=(1,0,\sqrt3)$. When $\alpha=2\sqrt2-2$, one can verify that
$${\rm card}\{X(\alpha,\Lambda)\}=20.$$

By Remark 2.1, for a packing lattice $\Lambda$ of $B^3$ and $\alpha=2\sqrt2-2$, to let
$${\rm card}\{X(\alpha,\Lambda)\}\geq20,$$
a necessary condition is there exist ${\bf a}_1, {\bf a}_2\in\Lambda$ such that
$$||{\bf a}_1||=||{\bf a}_2||=2,\quad ||{\bf a}_1+{\bf a}_2||=||{\bf a}_1-{\bf a}_2||=2\sqrt2.$$
Without loss of generality, we suppose ${\bf a}_1=(2,0,0),\ {\bf a}_2=(0,2,0)$. It is obvious to see that we can expand a basis of $\Lambda$ based on $\{{\bf a}_1, {\bf a}_2\}$.

Suppose ${\bf a}_1, {\bf a}_2$ and ${\bf a}_3=(v_1, v_2, v_3)$ is a basis of $\Lambda$. We assume
$$v_3>0,$$
without loss of generality. Since
$$\det(\Lambda)\geq\frac{{\rm vol}(B^3)}{\delta^*(B^3)}=\frac{\frac{4}{3}\pi}{\frac{\pi}{\sqrt{18}}}=4\sqrt2,$$
we have
$$v_3\geq\sqrt2$$
and the equality holds if and only if $\Lambda$ is the densest packing lattice of $B^3$. In this case one can verify that
$${\rm card}\{X(\alpha,\Lambda)\}=18<20.$$
Therefore, we have
$$v_3>\sqrt2.$$
Which means that, for
$${\bf v}=z_1{\bf a}_1+z_2{\bf a}_2+z_3{\bf a}_3\in X(\alpha,\Lambda),$$
we have
$$z_3=0\ {\rm or}\ \pm1.$$

Since
$${\rm card}\{{\bf v}:\ {\bf v}=z_1{\bf a}_1+z_2{\bf a}_2\in X(\alpha,\Lambda)\}=8$$
and $X(\alpha,\Lambda)$ is centrally symmetric, to let ${\rm card}\{X(\alpha,\Lambda)\}\geq20$, we have
$${\rm card}\{{\bf v}:\ {\bf v}=z_1{\bf a}_1+z_2{\bf a}_2+{\bf a}_3\in X(\alpha,\Lambda)\}\geq6.$$
Therefore, there exist two of them is equivalent. Replace ${\bf a}_3$ by the mid-point of them, we may further assume
$$||{\bf a}_3||^2=v_1^2+v_2^2+v_3^2=4,$$
by Remark 2.1.

Without loss of generality, we suppose $v_1, v_2\geq0$. Since
$$||{\bf a}_3-{\bf a}_1||\geq2,\quad ||{\bf a}_3-{\bf a}_2||\geq2,$$
by routine computation we have
$$0\leq v_1\leq1,\quad 0\leq v_2\leq1.$$

For a lattice vector
$${\bf v}=z_1{\bf a}_1+z_2{\bf a}_2+{\bf a}_3\in X(\alpha,\Lambda),$$
we have
$$||{\bf v}||^2=4z_1^2+4z_1v_1+4z_2^2+4z_2v_2+v_1^2+v_2^2+v_3^2\leq(2\sqrt2)^2,$$
which means
$$z_1^2+z_1v_1+z_2^2+z_2v_2\leq1. \eqno(4.9)$$

By routine computation, a necessary condition for (4.9) is
$$|z_1|\leq1,\quad |z_2|\leq1.$$
For $(z_1, z_2)=(1, 0)$: ${\bf a}_1+{\bf a}_3\in X(\alpha,\Lambda)$ if and only if
$$v_1=0. \eqno(4.9.1)$$
For $(z_1, z_2)=(1, -1)$: ${\bf a}_1-{\bf a}_2+{\bf a}_3\in X(\alpha,\Lambda)$ if and only if
$$1+v_1-v_2\leq0. \eqno(4.9.2)$$
For $(z_1, z_2)=(0, 1)$: ${\bf a}_2+{\bf a}_3\in X(\alpha,\Lambda)$ if and only if
$$v_2=0. \eqno(4.9.3)$$
For $(z_1, z_2)=(-1, 1)$: $-{\bf a}_1+{\bf a}_2+{\bf a}_3\in X(\alpha,\Lambda)$ if and only if
$$1-v_1+v_2\leq0. \eqno(4.9.4)$$
For $(z_1, z_2)=(-1, -1)$: $-{\bf a}_1-{\bf a}_2+{\bf a}_3\in X(\alpha,\Lambda)$ if and only if
$$1-v_1-v_2\leq0. \eqno(4.9.5)$$
For $(z_1, z_2)=(1, 1)$, since $1+v_1+1+v_2>1$, we have
$${\bf a}_1+{\bf a}_2+{\bf a}_3\notin X(\alpha,\Lambda). \eqno(4.9.6)$$

Obviously, (4.9.1) and (4.9.4) cannot hold simultaneously, (4.9.2) and (4.9.3) cannot hold simultaneously. Combining with (4.9.6), we have
$${\rm card}\{{\bf v}:\ {\bf v}=z_1{\bf a}_1+z_2{\bf a}_2+{\bf a}_3\in X(\alpha,\Lambda)\}\leq6.$$
Therefore, for $\alpha=2\sqrt2-2$ we get
$${\rm card}\{X(\alpha,\Lambda)\}\leq20.$$

To let ${\rm card}\{X(\alpha,\Lambda)\}=20$, we must have: one of (4.9.1) and (4.9.4) holds, one of (4.9.2) and (4.9.3) holds, and (4.9.5) holds. By routine computation one can deduce that
$${\bf a}_3=(0, 1, \sqrt3)\ \ {\rm or}\ \ (1, 0, \sqrt3).$$
Therefore, when $\alpha=2\sqrt2-2$, we have
$$\kappa^*_\alpha(B^3)=20,$$
and the equality holds if and only if the corresponding lattice $\Lambda$ is generated by ${\bf a}_1=(2, 0, 0)$, ${\bf a}_2=(0, 2, 0)$ and ${\bf a}_3=(1, 0, \sqrt3)$, up to rotation and reflection.
Lemma 4.5 is proved.
\end{proof}

\medskip
Lemma 4.1, Lemma 4.4 and Lemma 4.5 together yields the following theorem.

\medskip
\noindent
{\bf Theorem 2.} {\it In $\mathbb{E}^3$, we have}
$$\kappa^*_\alpha(B^3)=\left\{
\begin{array}{ll}
12,&  0\leq\alpha<\frac{4}{3}\sqrt3-2,\\
14,& \frac{4}{3}\sqrt3-2\leq\alpha<2\sqrt2-2,\\
20,& \alpha=2\sqrt2-2.
\end{array}
\right.
$$

\medskip
\noindent
{\bf Remark 4.2.} In fact, by repeating the calculations (4.9.1)-(4.9.6), one can deduce that for $\alpha=2\sqrt2-2$,
$${\rm card}\{X(\alpha,\Lambda)\}=18$$
if and only if $\Lambda$ is generated by ${\bf a}_1=(2, 0, 0)$, ${\bf a}_2=(0, 2, 0)$ and ${\bf a}_3=(1, 1, \sqrt2)$, or ${\bf a}_1=(2, 0, 0)$, ${\bf a}_2=(0, 2, 0)$ and ${\bf a}_3=(0, 0, 2)$, up to rotation and reflection.

\medskip
We end this section by a problem as following.

\medskip
\noindent
{\bf Problem 4.1.} {\it When $\alpha=\frac{4}{3}\sqrt3-2$, is it true that
$$\kappa^*_{\alpha}(C)\geq14$$
holds for every three-dimensional centrally symmetric convex body $C$?}

\bigskip
\noindent
{\LARGE\bf 5. Kissing numbers of convex bodies}

\bigskip
Although the concept of $\kappa^*_\alpha(B^3)$ itself is interesting, it can also leads to determine the lattice kissing numbers of convex bodies which were geometrically similar to $B^3$. To this end, we present the following theorem:

\medskip
\noindent
{\bf Theorem 5.1.} {\it For a 3-dimensional centrally symmetric convex body $C$ centered at ${\bf o}$, if $$B^3\subset C\subset {\rm int}(\mbox{$\frac{2}{3}$}\sqrt3B^3)$$
holds, then we have
$$\kappa^{*}(C)=12.$$}
\begin{proof}
Let $C+\Lambda$ be a lattice packing attaining $\kappa^{*}(C)$ and
$$X=\{{\bf v}_1,...,{\bf v}_{\kappa^*(C)}\}=\partial(2C)\cap\Lambda,$$
where $\partial(2C)$ denotes the boundary of $2C$. Since
$$2B^3\subset 2C\subset {\rm int}(\mbox{$\frac{4}{3}$}\sqrt3B^3),$$
we have
$$2\leq||{\bf v}_i||<\mbox{$\frac{4}{3}$}\sqrt3$$
holds for all $i=1,2,...,\kappa^*(C)$. Since $\Lambda$ is also a packing lattice of $B^3$, by Lemma 4.3 we get
$$\kappa^*(C)={\rm card}X\leq12.$$

On the other hand, since
$$\kappa^*(C)\geq12$$
holds for all 3-dimensional centrally symmetric convex body $C$ (see \cite{Sw01}), we have
$$\kappa^{*}(C)=12.$$

Theorem 5.1 is proved.
\end{proof}

\medskip
We now give several convex bodies as examples which lattice kissing numbers can be determined by Theorem 5.1.

\medskip
\noindent
{\bf Example 5.1.} We take $\tau=\frac{\sqrt5+1}{2}$ and define
$$P_d=\{(v_1,v_2,v_3): |\tau v_1|+|v_2|\leq1,\ |\tau v_2|+|v_3|\leq1,\ |\tau v_3|+|v_1|\leq1\},$$
\begin{eqnarray*}
P_i=\{(v_1,v_2,v_3)\hspace{-0.3cm}&:&\hspace{-0.2cm} |v_1|+|v_2|+|v_3|\leq1,\ |\tau v_1|+|\mbox{$\frac{1}{\tau}$}v_3|\leq1,\ |\tau v_2|+|\mbox{$\frac{1}{\tau}$}v_1|\leq1,\\
& &\hspace{-0.2cm}|\tau v_3|+|\mbox{$\frac{1}{\tau}$}v_2|\leq1\}.
\end{eqnarray*}
Usually, $P_d$ and $P_i$ are called a dodecahedron and an icosahedron, respectively. Define
$$P_{tri}=(1+\tau)P_i\cap(4/3+\tau)P_d.$$
Usually, $P_{tri}$ is called a truncated icosahedron.

\begin{figure}[h]
\centering
\includegraphics[height=6.5cm]{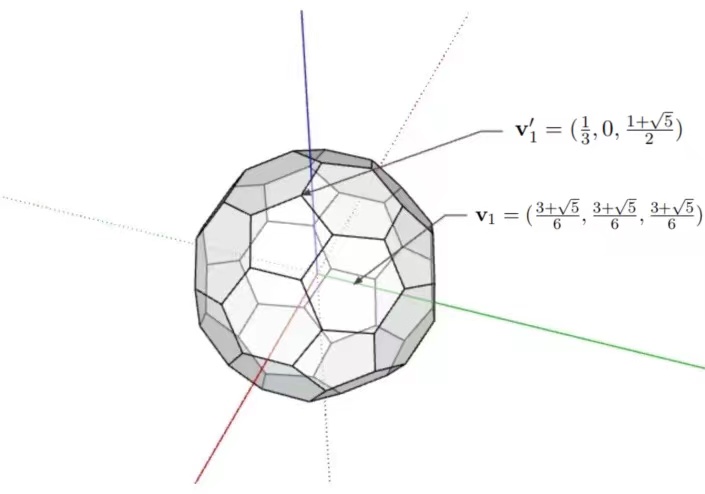}
\caption{Truncated icosahedron $P_{tri}$}
\label{fig1}
\end{figure}

By routine computation, we have
$$||{\bf v}_1||B^3\subset P_{tri}\subset||{\bf v}_1'||B^3,$$
where ${\bf v}_1=(\frac{3+\sqrt5}{6},\frac{3+\sqrt5}{6},\frac{3+\sqrt5}{6})$, ${\bf v}_1'=(\frac{1}{3},0,\frac{1+\sqrt5}{2})$, see Figure 1.
Since
$$B^3\subset\frac{1}{||{\bf v}_1||}P_{tri}\subset\frac{||{\bf v}_1'||}{||{\bf v}_1||}B^3=1.0929...B^3\subset {\rm int}(\mbox{$\frac{4}{3}$}\sqrt3B^3),$$
by Theorem 5.1 we have:
$$\kappa^*(P_{tri})=\kappa^*(\frac{1}{||{\bf v}_1||}P_{tri})=12.$$

\medskip
\noindent
{\bf Example 5.2.} We define
\begin{eqnarray*}
P_{rtc}=\{(v_1,v_2,v_3)\hspace{-0.3cm}&:&\hspace{-0.3cm}|\tau v_1|\leq1,|\tau v_2|\leq1,|\tau v_3|\leq1,|\frac{1}{2}v_1|+|\frac{\tau}{2}v_2|+|\frac{\tau+1}{2}v_3|\leq1,\\
& &|\frac{\tau}{2}v_1|+|\frac{\tau+1}{2}v_2|+|\frac{1}{2}v_3|\leq1,|\frac{\tau+1}{2}v_1|+|\frac{1}{2}v_2|+|\frac{\tau}{2}v_3|\leq1\}.
\end{eqnarray*}
Usually, $P_{rtc}$ is called a rhombic triacontahedron. Define
$$P_{rid}=(3\tau+2)P_{rtc}\cap(4\tau+1)P_i\cap(3(1+\tau))P_d,$$
$$P_{trid}=(5\tau+4)P_{rtc}\cap(6\tau+3)P_i\cap(5(1+\tau))P_d.$$
Usually, $P_{rid}$ and $P_{trid}$ are called a rhombic icosidodecahedron and a truncated icosidodecahedron, respectively.

\begin{figure}[h]
\raggedright
\includegraphics[height=4.6cm]{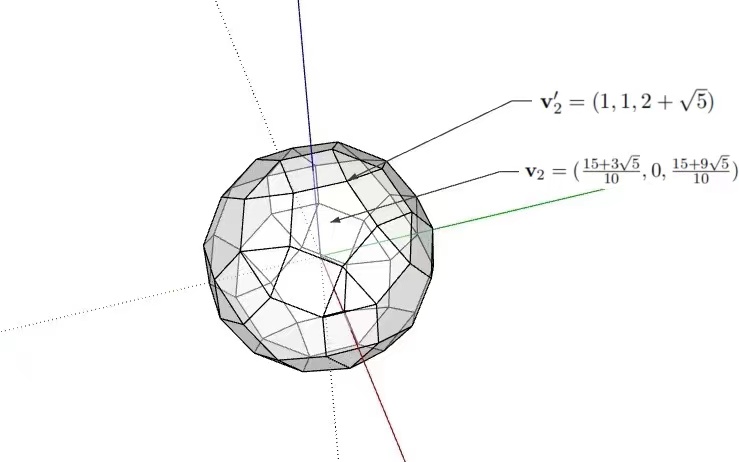}
\label{fig 2}
\raggedleft
\includegraphics[height=4.6cm]{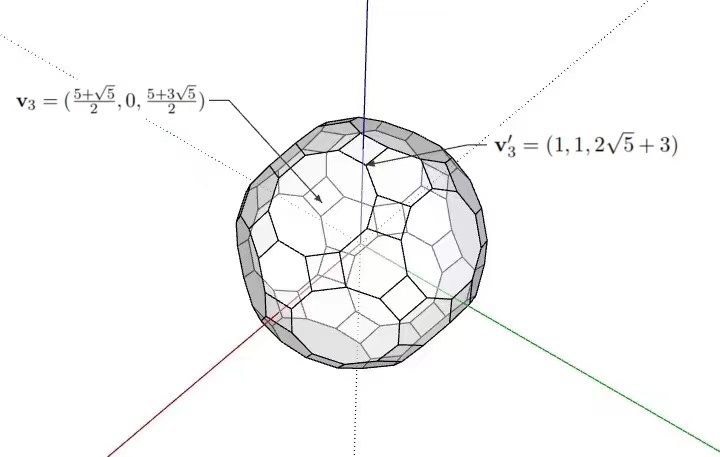}
\caption{Rhombic icosidodecahedron $P_{rid}$ and truncated icosidodecahedron $P_{trid}$}
\label{fig 2}
\end{figure}

By routine computation, we have
$$||{\bf v}_2||B^3\subset P_{rid}\subset||{\bf v}_2'||B^3,$$
where ${\bf v}_2=(\frac{15+3\sqrt5}{10},0,\frac{15+9\sqrt5}{10})$, ${\bf v}_2'=(1,1,2+\sqrt5)$, see Figure 2.
Since
$$B^3\subset\frac{1}{||{\bf v}_2||}P_{rid}\subset\frac{||{\bf v}_2'||}{||{\bf v}_2||}B^3=1.0815...B^3\subset {\rm int}(\mbox{$\frac{4}{3}$}\sqrt3B^3),$$
by Theorem 5.1 we have:
$$\kappa^*(P_{rid})=\kappa^*(\frac{1}{||{\bf v}_2||}P_{rid})=12.$$

By routine computation, we have
$$||{\bf v}_3||B^3\subset P_{trid}\subset||{\bf v}_3'||B^3,$$
where ${\bf v}_3=(\frac{5+\sqrt5}{2},0,\frac{5+3\sqrt5}{2})$, ${\bf v}_3'=(1,1,2\sqrt5+3)$, see Figure 2.
Since
$$B^3\subset\frac{1}{||{\bf v}_3||}P_{trid}\subset\frac{||{\bf v}_3'||}{||{\bf v}_3||}B^3=1.1050...B^3\subset {\rm int}(\mbox{$\frac{4}{3}$}\sqrt3B^3),$$
by Theorem 5.1 we have:
$$\kappa^*(P_{trid})=\kappa^*(\frac{1}{||{\bf v}_3||}P_{trid})=12.$$

\medskip
\noindent
{\bf Example 5.3.} We use the configuration of snub dodecahedron which given by M. Henk \cite{He01}, denote it by $P_{sd}$. By routine computation, we have
$$||{\bf v}_4||B^3\subset P_{sd}\subset||{\bf v}_4'||B^3,$$
where
$${\bf v}_4=(-0.9661..., 0, 1.5632...),\ {\bf v}_4'=(-0.3477..., -0.3069..., 1.9454...),$$
see Figure 3. Since
$$B^3\subset\frac{1}{||{\bf v}_4||}P_{sd}\subset\frac{||{\bf v}_4'||}{||{\bf v}_4||}B^3=1.0883...B^3\subset {\rm int}(\mbox{$\frac{4}{3}$}\sqrt3B^3),$$
by Theorem 5.1 we have:
$$\kappa^*(P_{sd})=\kappa^*(\frac{1}{||{\bf v}_4||}P_{sd})=12.$$

\begin{figure}[h]
\centering
\includegraphics[height=6.5cm]{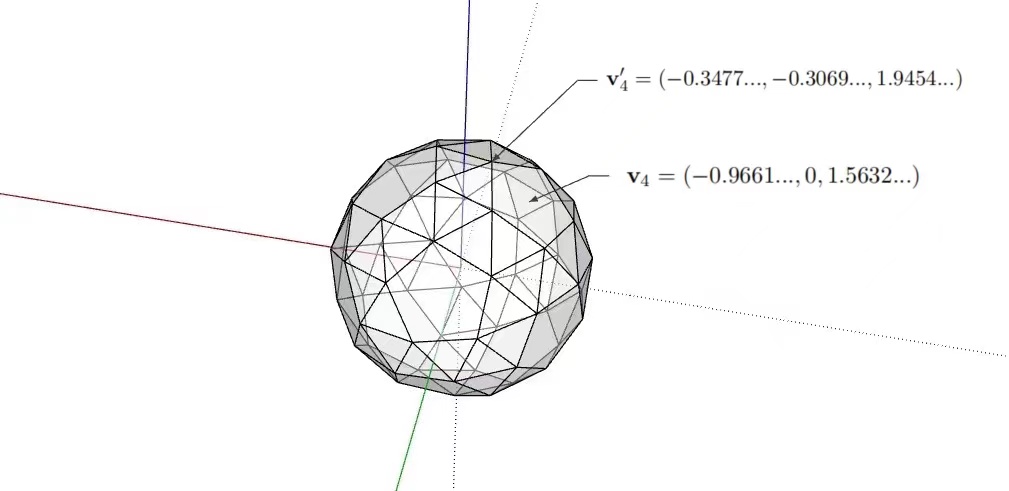}
\caption{Snub dodecahedron $P_{sd}$}
\label{fig 3}
\end{figure}

\medskip
\noindent
{\bf Example 5.4.} Define
$$B_p^{3}=\{(v_1,v_2,v_3):\ |v_1|^p+|v_2|^p+|v_3|^p\leq1\}.$$
Usually, $B_p^{3}$ is called a $L_p$ unit ball in 3-dimension. When $p_1\leq p_2$, it is well known that
$$B_{p_1}^{3}\subset B_{p_2}^{3}.$$

For $p\geq2$ and a point ${\bf v}=(v_1,v_2,v_3)\in B_p^{3}$ where $v_1,v_2,v_3\geq0$, we have
$$v_1^p+v_2^p+v_3^p\leq1.$$
According to Power-Mean Inequality, we have
$$v_1^2+v_2^2+v_3^2\leq{((\frac{1}{3})^{\frac{1}{p}})}^{2}\times3=(\frac{1}{3})^{\frac{2}{p}}\times3.$$
Therefore we have
$$B^3\subset B_p^{3}\subset\sqrt3\times(\frac{1}{3})^{\frac{1}{p}}B^3.$$
By routine computation, when $2\leq p<\frac{\ln3}{\ln3-\ln2}$ we get
$$\kappa^*(B_p^{3})=12$$
by Theorem 5.1.

For $p<2$, by the same deduction we obtain
$$\sqrt3\times(\frac{1}{3})^{\frac{1}{p}}B^3\subset B_p^{3}\subset B^3,$$
which means
$$B^3\subset \frac{1}{\sqrt3\times(\frac{1}{3})^{\frac{1}{p}}}B_p^{3}\subset \frac{1}{\sqrt3\times(\frac{1}{3})^{\frac{1}{p}}}B^3.$$
By routine computation, when $2> p>\frac{\ln3}{\ln2}$ we get
$$\kappa^*(B_p^{3})=12$$
by Theorem 5.1.

As a conclusion of the two cases, when $\frac{\ln3}{\ln3-\ln2}> p>\frac{\ln3}{\ln2}$, we have
$$\kappa^*(B_p^{3})=12.$$

\medskip
\noindent
{\bf Remark 5.1.} Let $\Lambda$ be the lattice generated by ${\bf a}_1=(2,0,0), {\bf a}_2=(0,2,0)$ and ${\bf a}_3=(1,1,1)$, then one can verify that: when $p=\frac{\ln3}{\ln2}$, $\Lambda$ is a packing lattice of $B_p^{3}$ and
$${\rm card}\{\partial(2B_p^{3})\cap\Lambda\}=14.$$
On the other hand, since $B_p^3$ is a strictly convex body when $1<p<\infty$, combining with Lemma 2.2, when $p=\frac{\ln3}{\ln2}$ we have
$$\kappa^*(B_p^{3})=14.$$

\bigskip
\noindent
{\LARGE\bf 6. Proof of Theorem 4}

\bigskip
For $\alpha=2\sqrt2-2$ and a packing lattice $\Lambda$ of $B^n$, we have the following lemma.

\medskip
\noindent
{\bf Lemma 6.1.} One equivalent class of $\Lambda$ can contain at most $n$ pairs of vectors of $X(\alpha,\Lambda)$.
\begin{proof}
Suppose
$$\pm{\bf v}_1,\ldots,\pm{\bf v}_i \in X(\alpha,\Lambda)$$
belong to the same equivalent class, $i\geq2$. By Remark 2.1, we have
$$||{\bf v}_i||=2\sqrt2$$
holds for all $i$ and
$$\langle{\bf v}_i, {\bf v}_j\rangle=0$$
holds for all $i\neq j$. Therefore, one equivalent class of $\Lambda$ can contain at most $n$ pairs of vectors of $X(\alpha,\Lambda)$.

Lemma 6.1 is proved.
\end{proof}

\medskip
Denote the numbers of equivalent classes of $\Lambda$ which contain exactly $i$ pairs of vectors of $X(\alpha,\Lambda)$ by $m_i$. We define a collection of sets
$$C(X(\alpha,\Lambda))=\big\{ A=\{{\bf v}_1, {\bf v}_2, {\bf v}_3\}:\ {\bf v}_2=\frac{1}{2}({\bf v}_1+{\bf v}_3)\ {\rm and} \ {\bf v}_1,\ {\bf v}_2,\ {\bf v}_3\in X(\alpha,\Lambda)\big\}.$$
By estimate ${\rm card}\{C(X(\alpha,\Lambda))\}$ in two different ways, we prove the following lemma.

\medskip
\noindent
{\bf Lemma 6.2.}
$$\sum\limits_{i=2}^{n}2i(i-1)m_i\leq \kappa^*(B^{n-1})\cdot m_1.$$

\begin{proof}
For a set $\{{\bf v}_1, {\bf v}_2, {\bf v}_3\}\in C(X(\alpha,\Lambda))$, by Remark 2.1, we have
$$||{\bf v}_2||=2,\ ||{\bf v}_1-{\bf v}_2||=||{\bf v}_3-{\bf v}_2||=2,\ ||{\bf v}_1||=||{\bf v}_3||=2\sqrt2.$$
We assume ${\bf v}_2=(0,0,\ldots,0,2)$, without loss of generality. Then one can easily deduce that ${\bf v}_1$ and ${\bf v}_3$ must lie in the $(n-1)$-dimensional hyperplane
$$\pi_0:\ \big\{(v_1,v_2,\ldots,v_{n-1},v_n):\ v_n=2\big\}.$$

It is obvious that
$${\rm card}\{{\bf v}: ||{\bf v}-{\bf v}_2||=2,\ {\bf v}\in X(\alpha,\Lambda)\cap\pi_0\}\leq\kappa^*(B^{n-1}),$$
which means
$${\rm card}\{A\in C(X(\alpha,\Lambda)):\ {\bf v}_2\in A\}\leq\kappa^*(B^{n-1})/2.$$

For a vector ${\bf v}\in X(\alpha,\Lambda)$ which length is $2$, by Remark 2.1, ${\bf v}$ cannot equivalent with vector of $X(\alpha,\Lambda)$, besides $\pm{\bf v}$. Therefore, we have
$${\rm card}\{{\bf v}: ||{\bf v}||=2,\ {\bf v}\in X(\alpha,\Lambda)\}\leq2m_1.$$
Consequently, we get
$${\rm card}\{C(X(\alpha,\Lambda))\}\leq\kappa^*(B^{n-1})\cdot m_1.\eqno(6.1)$$

On the other hand, by the definition of $C(X(\alpha,\Lambda))$, a set $\{{\bf v}_1, {\bf v}_2, {\bf v}_3\}$ belong to it if and only if ${\bf v}_1, {\bf v}_3\in X(\alpha,\Lambda)$ is equivalent and
$${\bf v}_1\neq\pm{\bf v}_3.$$
For an equivalent class which contain $i\geq2$ pairs of vectors of $X(\alpha,\Lambda)$, denote it by $X_{i1}$. By enumerate we have
$${\rm card}\big\{\{{\bf v}_1,{\bf v}_3\}: \{{\bf v}_1,{\bf v}_3\}\subset X_{i1},\ {\bf v}_1\neq\pm{\bf v}_3\big\}=2i(i-1).$$ Therefore, we get
$${\rm card}\{C(X(\alpha,\Lambda))\}=\sum\limits_{i=2}^{n}2i(i-1)m_i.\eqno(6.2)$$

By (6.1) and (6.2), Lemma 6.2 is proved.
\end{proof}

\medskip
\noindent
{\bf Theorem 4.} {\it In $\mathbb{E}^8$, when $\alpha=2\sqrt2-2$ we have}
$$\kappa^*_\alpha(B^8)=2400.$$

\begin{proof}
For $\alpha=2\sqrt2-2$, it is well known (see \cite{Co01}) that
$${\rm card}\{X(\alpha,\sqrt2E_8)\}=240+2160=2400,$$
where
$$E_8=\big\{(v_1,v_2,\ldots,v_8):\ 2v_i\in Z;\ v_i-v_j\in Z;\ \sum v_i\in 2Z\big\}.$$

Suppose that there are a suitable lattice $\Lambda$ satisfying
$${\rm card}\{X(\alpha,\Lambda)\}\geq2400,$$
which means
$$8m_8+7m_7+\ldots+2m_2+m_1\geq1200.\eqno(6.3)$$
Since there are at most $2^8-1=255$ equivalent classes can contain the vectors of $X(\alpha,\Lambda)$, we have
$$m_8+m_7+\ldots+m_1\leq255.\eqno(6.4)$$
For $n=8$ case, we restate Lemma 6.2 as
$$112m_8+84m_7+60m_6+40m_5+24m_4+12m_3+4m_2\leq126m_1\eqno(6.5)$$
by substitute $\kappa^*(B^7)=126$.

By (6.3) and (6.4), we have
$$7m_8+6m_7+\ldots+m_2\geq945,$$
multiply both sides by $34$, we have
$$238m_8+204m_7+170m_6+136m_5+102m_4+68m_3+34m_2\geq32130.\eqno(6.6)$$
By (6.4) and (6.5), one can deduce that
$$238m_8+210m_7+186m_6+166m_5+150m_4+138m_3+130m_2\leq32130.\eqno(6.7)$$
By (6.6) and (6.7), we obtain
$$m_7=m_6=\ldots=m_2=0.$$
Combining with (6.3), (6.4) and (6.5), we have
$$\left\{
\begin{aligned}
& 8m_8+m_1\geq1200,\\
& m_8+m_1\leq255,\\
& 112m_8\leq126m_1.\\
\end{aligned}
\right.
$$

By routine computation, one can easily deduce that
$$m_1=120,\ m_8=135.$$
Furthermore, in this case the equality in (6.1) holds, which means that:
$${\rm card}\{{\bf v}:||{\bf v}||=2,\ {\bf v}\in X(\alpha,\Lambda)\}=2m_1=240.$$
Since
$$\kappa^*(B^8)=240$$
and the corresponding lattice must be $\sqrt2E_8$, up to rotation and reflection (see \cite{Zo06}), for $\alpha=2\sqrt2-2$ we have
$$\kappa^*_\alpha(B^8)=2400,$$
and the equality can be attained if and only if the corresponding lattice $\Lambda=\sqrt2E_8$, up to rotation and reflection.

Theorem 4 is proved.
\end{proof}

\medskip
Based on this proof, we may make the following conjecture.

\medskip
\noindent
{\bf Conjecture 6.1.} {\it In $\mathbb{E}^8$, when $\alpha=2\sqrt2-2$ we have}
$$\kappa_\alpha(B^8)=2400.$$

\medskip
Let $\Lambda_{24}$ denote the {\it Leech lattice} (see \cite{Co01}). When $\alpha=2\sqrt2-2$, we have
$${\rm card}\{X(\alpha,\sqrt2\Lambda_{24})\}=196560+16773120+398034000=415003680.$$
This observation supports the following conjecture.

\medskip
\noindent
{\bf Conjecture 6.2.} {\it In $\mathbb{E}^{24}$, when $\alpha=2\sqrt2-2$ we have}
$$\kappa^*_\alpha(B^{24})=415003680.$$

\bigskip
\noindent
{\LARGE\bf 7. Proof of Theorem 3}

\bigskip
\noindent
{\bf Theorem 3.} {\it In $\mathbb{E}^4$, we have}
$$\kappa^*_\alpha(B^4)=\left\{
\begin{aligned}
30,& \ \ \sqrt6-2\leq\alpha<2\sqrt2-2,\\
50,& \ \ \alpha=2\sqrt2-2.\\
\end{aligned}
\right.
$$
\begin{proof}
As usual, we write
$$A_n=\big\{(v_1,v_2,\ldots,v_n):\ v_i\in Z;\ \sum v_i=0\big\},$$
$$D_n=\big\{(v_1,v_2,\ldots,v_n):\ v_i\in Z;\ \sum v_i\in 2Z\big\}.$$
Furthermore, we denote the {\it dual lattice} of $A_n$ by $A_n^*$, namely
$$A_n^*=\big\{{\bf v}:\ \langle{\bf v}, {\bf u}\rangle\in Z\ for\ all\ {\bf u}\in A_n\big\}.$$

When $\alpha=\sqrt6-2$, one can verify that
$${\rm card}\{X(\alpha,\sqrt5A_4^*)\}=30.$$
Combining with Lemma 2.2, for $\sqrt6-2\leq\alpha<2\sqrt2-2$ we have
$$\kappa^*_\alpha(B^4)=30.$$

For $\alpha=2\sqrt2-2$, let $\Lambda$ be the lattice generated by ${\bf a}_1=(2,0,0,0)$, ${\bf a}_2=(0,2,0,0)$, ${\bf a}_3=(1,0,\sqrt3,0)$, ${\bf a}_4=(0,1,\frac{2}{3}\sqrt3, \frac{\sqrt5}{\sqrt3})$. One can verify that
$${\rm card}\{X(\alpha,\Lambda)\}=50.\eqno(7.1)$$

Suppose that there exist a packing lattice $\Lambda$ of $B^4$ satisfying
$${\rm card}\{X(\alpha,\Lambda)\}\geq52.$$
We still denotes the numbers of equivalent classes of $\Lambda$ which contain exactly $i$ pairs of vectors of $X(\alpha,\Lambda)$ by $m_i$.

If $m_4\neq0$, by Remark 2.1, we may assume
$${\bf v}_1=(2\sqrt2,0,0,0),\quad {\bf v}_2=(0,2\sqrt2,0,0),\quad {\bf v}_3=(0,0,2\sqrt2,0),\quad {\bf v}_4=(0,0,0,2\sqrt2)$$
belong to $X(\alpha,\Lambda)$ and
$$\frac{1}{2}({\bf v}_i\pm{\bf v}_j),\quad i\neq j$$
belong to $X(\alpha,\Lambda)$, without loss of generality. In this case, lattice $\Lambda$ is generated by
$${\bf a}_1=(\sqrt2,\sqrt2,0,0),\ {\bf a}_2=(\sqrt2,-\sqrt2,0,0),\ {\bf a}_3=(0,0,\sqrt2,\sqrt2),\ {\bf a}_4=(0,0,\sqrt2, -\sqrt2),$$
which means $\Lambda=\sqrt2D_4$. One can verify that
$${\rm card}\{X(\alpha,\sqrt2D_4)\}=24+24=48<52,$$
therefore
$$m_4=0.$$
Since $\sqrt2D_4$ is the unique densest packing lattice for $B^4$, up to rotation and reflection (see \cite{Zo06}), from now on we suppose
$$\det(\Lambda)>8.\eqno(7.2)$$

If $m_3=0$ and for every vector ${\bf v}\in X(\alpha,\Lambda)$ which length is $2$ we have
$${\rm card}\{A\in C(X(\alpha,\Lambda)):\ {\bf v}\in A\}<\kappa^*(B^3)/2=6,$$
by restate (6.1), (6.2), (6.3) and (6.4) for $n=4$, we obtain
$$\left\{
\begin{aligned}
& 4m_2\leq10m_1,\\
& 2m_2+m_1\geq26,\\
& m_2+m_1\leq15,\\
\end{aligned}
\right.
$$
which do not exist solution. Therefore, we have $m_3\neq0$ or there exist a vector ${\bf v}\in X(\alpha,\Lambda)$ which length is $2$ satisfy
$${\rm card}\{A\in C(X(\alpha,\Lambda)):\ {\bf v}\in A\}=\kappa^*(B^3)/2=6.$$

If $m_3\neq0$, by Remark 2.1, we assume
$${\bf v}_1=(2,-2,0,0),\quad {\bf v}_2=(2,2,0,0),\quad {\bf v}_3=(0,0,2\sqrt2,0)$$
belong to $X(\alpha,\Lambda)$ and
$$\frac{1}{2}({\bf v}_i\pm{\bf v}_j),\quad i\neq j$$
belong to $X(\alpha,\Lambda)$, without loss of generality. Therefore, the basis of lattice $\Lambda$ can be expanded by
$${\bf a}_1=(2,0,0,0),\quad {\bf a}_2=(0,2,0,0),\quad {\bf a}_3=(1,1,\sqrt2,0).$$

On the other hand, if there exist a vector ${\bf v}\in X(\alpha,\Lambda)$ which length is $2$ satisfy
$${\rm card}\{A\in C(X(\alpha,\Lambda)):\ {\bf v}\in A\}=\kappa^*(B^3)/2=6,$$
then there exist a three-dimensional subspace $H_0$ satisfy
$${\rm card}\{{\bf v}: {\bf v}\in H_0\cap\Lambda,\ ||{\bf v}||=2\}=\kappa^*(B^3)=12.$$
Therefore, we may suppose
$$H_0=\{(v_1,v_2,v_3,v_4):\ v_4=0\}$$
and the three-dimensional lattice $H_0\cap\Lambda$ is generated by
$${\bf a}_1=(2,0,0,0),\quad {\bf a}_2=(0,2,0,0),\quad {\bf a}_3=(1,1,\sqrt2,0),$$
without loss of generality.

As a conclusion of two cases above, we set a basis of lattice $\Lambda$ by
$${\bf a}_1=(2,0,0,0),\ {\bf a}_2=(0,2,0,0),\ {\bf a}_3=(1,1,\sqrt2,0),\ {\bf a}_4=(v_1,v_2,v_3,v_4)$$
and
$$v_1\geq0,\quad v_2\geq0,\quad v_3\geq0,\quad v_4\geq0$$
without loss of generality. Furthermore, by (7.2) we have
$$v_4>\sqrt2.$$
Therefore, for a vector
$${\bf v}=z_1{\bf a}_1+z_2{\bf a}_2+z_3{\bf a}_3+z_4{\bf a}_4\in X(\alpha,\Lambda),$$
we have
$$z_4=0\ {\rm or}\ \pm1.$$

By Remark 4.2,
$${\rm card}\{{\bf v}:\ {\bf v}=z_1{\bf a}_1+z_2{\bf a}_2+z_3{\bf a}_3\in X(\alpha,\Lambda)\}=18.$$
Since $X(\alpha,\Lambda)$ is centrally symmetric, we have
$${\rm card}\{{\bf v}:\ {\bf v}=z_1{\bf a}_1+z_2{\bf a}_2+z_3{\bf a}_3+{\bf a}_4\in X(\alpha,\Lambda)\}\geq17,$$
which means that there exist two of them is equivalent. Replace ${\bf a}_4$ by the mid-point of them, we may further assume
$$||{\bf a}_4||^2=v_1^2+v_2^2+v_3^2+v_4^2=4,$$
by Remark 2.1.

By routine computation, besides $z_1{\bf a}_1+z_2{\bf a}_2+z_3{\bf a}_3$, vector which belong to $X(\alpha,\Lambda)$ must be one of the following form:
$$z_1{\bf a}_1+z_2{\bf a}_2\pm{\bf a}_4,\quad z_1{\bf a}_1+z_2{\bf a}_2\pm({\bf a}_3+{\bf a}_4),$$
$$z_1{\bf a}_1+z_2{\bf a}_2\pm({\bf a}_3-{\bf a}_4),\quad z_1{\bf a}_1+z_2{\bf a}_2\pm(2{\bf a}_3-{\bf a}_4).$$
To let ${\rm card}\{X(\alpha,\Lambda)\}\geq52$, there exist one form above have at least ten vectors which belongs to $X(\alpha,\Lambda)$. Without loss of generality, we suppose
$${\rm card}\{{\bf v}:{\bf v}=z_1{\bf a}_1+z_2{\bf a}_2\pm{\bf a}_4\in X(\alpha,\Lambda)\}\geq10.$$
Combining with
$${\rm card}\{{\bf v}:{\bf v}=z_1{\bf a}_1+z_2{\bf a}_2\in X(\alpha,\Lambda)\}=8,$$
by Lemma 4.4 and Remark 4.2, we may assume
$$v_1=0,\quad v_2=0,\quad v_3^2+v_4^2=4$$
or
$$v_1=0,\quad v_2=1,\quad v_3^2+v_4^2=3$$
without loss of generality.

For case $v_1=0,\ v_2=0,\ v_3^2+v_4^2=4$, by routine computation we have:
$${\rm card}\{{\bf v}:{\bf v}=z_1{\bf a}_1+z_2{\bf a}_2\pm{\bf a}_4\in X(\alpha,\Lambda)\}=10,$$
$${\rm card}\{{\bf v}:{\bf v}=z_1{\bf a}_1+z_2{\bf a}_2\pm ({\bf a}_3+{\bf a}_4)\in X(\alpha,\Lambda)\}=\left\{
\begin{aligned}
8,& \ \ v_3=0,\\
0,& \ \ v_3\neq0,\\
\end{aligned}
\right.
$$
$${\rm card}\{{\bf v}:{\bf v}=z_1{\bf a}_1+z_2{\bf a}_2\pm ({\bf a}_3-{\bf a}_4)\in X(\alpha,\Lambda)\}=8,$$
$${\rm card}\{{\bf v}:{\bf v}=z_1{\bf a}_1+z_2{\bf a}_2\pm (2{\bf a}_3-{\bf a}_4)\in X(\alpha,\Lambda)\}=\left\{
\begin{aligned}
2,& \ \ v_3\geq1/\sqrt2,\\
0,& \ \ v_3=0.\\
\end{aligned}
\right.
$$
Therefore, in this case we have
$${\rm card}\{X(\alpha,\Lambda)\}\leq44.$$

For case $v_1=0,\ v_2=1,\ v_3^2+v_4^2=3$, since
$$||{\bf a}_3-{\bf a}_4||\geq2,$$
we have
$$v_3\leq1/\sqrt2.$$
By routine computation we have:
$${\rm card}\{{\bf v}:{\bf v}=z_1{\bf a}_1+z_2{\bf a}_2\pm{\bf a}_4\in X(\alpha,\Lambda)\}=12,$$
$${\rm card}\{{\bf v}:{\bf v}=z_1{\bf a}_1+z_2{\bf a}_2\pm ({\bf a}_3+{\bf a}_4)\in X(\alpha,\Lambda)\}=4,$$
$${\rm card}\{{\bf v}:{\bf v}=z_1{\bf a}_1+z_2{\bf a}_2\pm ({\bf a}_3-{\bf a}_4)\in X(\alpha,\Lambda)\}=\left\{
\begin{aligned}
12,& \ \ v_3=1/\sqrt2,\\
4,& \ \ v_3<1/\sqrt2,\\
\end{aligned}
\right.
$$
$${\rm card}\{{\bf v}:{\bf v}=z_1{\bf a}_1+z_2{\bf a}_2\pm (2{\bf a}_3-{\bf a}_4)\in X(\alpha,\Lambda)\}=\left\{
\begin{aligned}
4,& \ \ v_3=1/\sqrt2,\\
0,& \ \ v_3<1/\sqrt2.\\
\end{aligned}
\right.
$$
Therefore, in this case we have ${\rm card}\{X(\alpha,\Lambda)\}\leq50$.

As a conclusion of these two cases and (7.1), for $\alpha=2\sqrt2-2$ we have
$$\kappa^*_\alpha(B^4)=50.$$
Theorem 3 is proved.
\end{proof}

\medskip
\noindent
{\bf Remark 7.1.} It is interesting to see that the $\sqrt2D_4$ lattice is not the optimal lattice in this case. Let $\Lambda$ be the lattice generated by ${\bf a}_1=(2,0,0,0)$, ${\bf a}_2=(0,2,0,0)$, ${\bf a}_3=(1,0,\sqrt3,0)$ and ${\bf a}_4=(0,1,0,\sqrt3)$. It is easy to show that, when $\alpha=2\sqrt2-2$, 
$${\rm card}\{X(\alpha,\Lambda)\}={\rm card}\{X(\alpha,\sqrt2D_4)\}=48.$$

\bigskip
\noindent
{\LARGE\bf 8. A Link Between $\kappa^*_\alpha(B^n)$ and $\gamma^*(B^n)$}

\bigskip
In 1964, P. Erd\"{o}s and C. A. Rogers \cite{Er01} studied the star number of the lattice covering for a convex body and proved the following result.

\medskip
\noindent
{\bf Theorem 8.1.} {\it Let $C$ be an ${\bf o}$-symmetric strictly convex body and $\Lambda$ a covering lattice of $C$ in $\mathbb{E}^n$. Then the star number of the covering $\{C+{\bf v}:\ {\bf v}\in  \Lambda\}$ is at least $2^{n+1}-1$, where the star number is the numbers of the translates of $C$ by lattice vectors, including $C$, which intersect the body $C$.}

\medskip
Let $\gamma^*(B^n)$ be the lattice packing-covering constant of $B^n$, namely
$$\gamma^{*}(B^n)=\min_\Lambda\{r: rB^n+\Lambda\ {\rm is \ a \ covering \ of} \ {\mathbb E^n}\},$$
where $\Lambda$ is a lattice such that $B^n+\Lambda$ is a packing in ${\mathbb E^n}$. For more details about $\gamma^*(B^n)$, we refer to \cite{Zo04}.

There exist a strong relation between $\kappa_\alpha^*(B^n)$ and $\gamma^*(B^n)$:

\medskip
\noindent
{\bf Theorem 8.2.} {\it For a given dimension $n_0$, suppose $\gamma^*(B^{n_0})=\sqrt2-\beta$ for a positive number $\beta$. Then for $\alpha\in[2\sqrt2-2\beta-2,\ 2\sqrt2-2)$ we have
$$\kappa_\alpha^*(B^{n_0})=2^{n_0+1}-2.$$
Which means that, if
$$\kappa_\alpha^*(B^{n_0})<2^{n_0+1}-2$$
holds for $\alpha<2\sqrt2-2$, then we have
$$\gamma^*(B^{n_0})\geq\sqrt2.$$}

\begin{proof}
We assume that $B^{n_0}+\Lambda$ is a lattice packing attaining $\gamma^*(B^{n_0})=\sqrt2-\beta$ for a positive $\beta$. For convenience, let
$$X=\{{\bf v}:\ 2\leq||{\bf v}||\leq2\sqrt2-2\beta,\ {\bf v}\in\Lambda\}.$$

It is easy to see that the star number of the covering configuration $(\sqrt2-\beta)B^{n_0}+\Lambda$ is ${\rm card}X+1$. By Theorem 8.1 we have
$${\rm card}X\geq 2^{n_0+1}-2.$$
Combining with Lemma 2.2, we get
$${\rm card}X=2^{n_0+1}-2.$$

Therefore, for $\alpha\in[2\sqrt2-2\beta-2,\ 2\sqrt2-2)$, we have
$$\kappa_\alpha^*(B^{n_0})=2^{n_0+1}-2.$$
Theorem 8.2 is proved.
\end{proof}

\medskip
\noindent
{\bf Remark 8.1.} Notice that
$$\gamma^*(B^5)>\sqrt2,$$
see \cite{Zo04}. However, when $\alpha=2\sqrt{9/5}-2$, one can verify that
$${\rm card}\{X(\alpha,\sqrt{24/5}A_5^*)\}=62.$$
Combining with Lemma 2.2, when $\alpha\in[2\sqrt{9/5}-2,\ 2\sqrt2-2)$ we have
$$\kappa_\alpha^*(B^5)=62.$$

\medskip
\noindent
{\bf Corollary 8.1.} In \cite{Sc02}, A. Sch\"{u}rmann and F. Vallentin improved the former best known result (see \cite{Zo04})
$$\gamma^*(B^6)\leq\sqrt2$$
to
$$\gamma^*(B^6)\leq2\sqrt{2\sqrt{798}-56}=1.411081242...$$
Therefore, by Theorem 8.2, for $\alpha\in[0.8222, 2\sqrt2-2)$ we have
$$\kappa_\alpha^*(B^6)=126.$$

\medskip
For a packing lattice $\Lambda$ of $B^n$ and an $\alpha<2\sqrt2-2$, the sufficient and necessary condition for
$${\rm card}\{X(\alpha,\Lambda)\}=2^{n+1}-2$$
is each equivalent class of $\Lambda$ must contain a pair of vectors of $X(\alpha,\Lambda)$. It is reasonable to imagine that, this condition is hard to satisfy in high dimensions. If so, the following conjecture make sense.

\medskip
\noindent
{\bf Conjecture 8.1.} {\it There are infinity numbers of dimension $n$ such that, when $\alpha<2\sqrt2-2$ we have
$$\kappa_\alpha^*(B^n)<2^{n+1}-2.$$
Especially, when $\alpha<2\sqrt2-2$ we have
$$\kappa_\alpha^*(B^8)<510$$
and
$$\kappa_\alpha^*(B^{24})<33554430.$$}

\medskip
\noindent
{\bf Remark 8.2.}
If Conjecture 8.1 is true, by Theorem 8.2, we have
$$\gamma^*(B^8)=\sqrt2$$
and
$$\gamma^*(B^{24})=\sqrt2,$$
which give an affirmative answer for C. Zong's Conjecture 3.1 in \cite{Zo04}.

\medskip
We write
$$E_7=\big\{{\bf v}:\ {\bf v}\in E_8;\ \sum v_i=0\big\}$$
and
$$E_6=\big\{{\bf v}:\ {\bf v}\in E_8;\ \sum v_i=v_7+v_8=0\big\}.$$
When $\alpha=2\sqrt2-2$, we have
$${\rm card}\{X(\alpha,\sqrt2D_5)\}=130,$$
$${\rm card}\{X(\alpha,\sqrt2E_6)\}=342,$$
and
$${\rm card}\{X(\alpha,\sqrt2E_7)\}=882,$$
see \cite{Co01}.

To end this article, we list some knowned results of $\kappa^*_\alpha(B^n)$ as the following table.

\medskip
\centerline{
\renewcommand\arraystretch{1.5}
\begin{tabular}{|c|c|c|}
\hline $n$ & $\kappa_\alpha^*(B^n)\ {\rm for}\ \alpha<2\sqrt2-2$ & $\kappa_\alpha^*(B^n)\ {\rm for}\ \alpha=2\sqrt2-2$\\
\hline $2$ & $\leq6 \ ({\rm can\ be\ attained})$ & $=8$\\
\hline $3$ & $\leq14 \ ({\rm can\ be\ attained})$ & $=20$\\
\hline $4$ & $\leq30 \ ({\rm can\ be\ attained})$ & $=50$\\
\hline $5$ & $\leq62 \ ({\rm can\ be\ attained})$ & $\geq130$\\
\hline $6$ & $\leq126 \ ({\rm can\ be\ attained})$ & $\geq342$\\
\hline $7$ & $\leq254 \ (??)$ & $\geq882$\\
\hline $8$ & $\leq510 \ (??)$ & $=2400$\\
\hline $24$ & $\leq33554430 \ (??)$ & $\geq 415003680$\\
\hline
\end{tabular}}

\vspace{0.8cm}\noindent
{\bf Acknowledgement.} This work is supported by the National Natural Science Foundation of China (NSFC12226006, NSFC11921001) and the Natural Key Research and Development Program of China (2018YFA0704701).

\vspace{0.6cm}
\noindent
Yiming Li, Center for Applied Mathematics, Tianjin University, Tianjin 300072, China
\noindent
Email: xiaozhuang@tju.edu.cn

\medskip\noindent
Chuanming Zong, Center for Applied Mathematics, Tianjin University, Tianjin 300072, China

\noindent
Email: cmzong@tju.edu.cn


\begin{thebibliography}{1}
\bibitem{Ba01}E. Bannai, N. J. A. Sloane, Uniqueness of certain spherical codes, \emph{Canadian J. Math.} \textbf{33} (1981), 437-449.
\bibitem{Co01}J. H. Conway, N. J. A. Sloane, \emph{Sphere packings, lattices and groups}, Springer-Verlag, New York, 1988.
\bibitem{Er01}P. Erd\"{o}s, C. A. Rogers, The star number of coverings of space with convex bodies, \emph{Acta Arith.} \textbf{9} (1964), 41-45.
\bibitem{Gu01}S. G\"{u}nter, Ein stereometrisches Problem, \emph{Grunert Arch.} \textbf{57} (1875), 209-215.
\bibitem{Ha01}H. Hadwiger, \"Uber Treffanzahlen bei translationgleichen Eik\"orpern, {\it Arch. Math.} {\bf 8} (1957), 212-213.
\bibitem{He01}M. Henk, http://www.eg-models.de/models/Lattices\_and\_Packings/Densest\_Polytopes/2001.02.060/\_applet.html, Electronic Geometry Model No. 2001.02.060.
\bibitem{Ho01}R. Hoppe, Bestimmung der gr\"{o}ssten Anzahl gleich grosser Kugeln, whlche sich auf eine Kugel von demselben Radius, wie die \"{u}brigen, auflegen lassen, \emph{Grunert Arch.} \textbf{56} (1874), 302-313.
\bibitem{Le02}J. Leech, The problem of the thirteen spheres, \emph{Math. Gazette} \textbf{40} (1956), 22-23.
\bibitem{Le01}V. I. Levenshtein, Boundaries for packings in $n$-dimensional Euclidean space, \emph{Dokl. Akad. Nauk SSSR} \textbf{245} (1979), 1299-1303.
\bibitem{Mi01}H. Minkowski, Dichteste gitterf\"ormige Lagerung kongruenter K\"orper, {\it Nachr. K. Ges. Wiss. G\"ottingen,} (1904), 311-355.
\bibitem{Mi02}H. Minkowski, \emph{Diophantische Approximationen}, Teubner, Leipzig 1907, 1927; Chelsea, New York, 1957.
\bibitem{Mu01}O. R. Musin, The problem of twenty-five spheres, \emph{Russian Math. Surveys} \textbf{58} (2003), 794-795.
\bibitem{Mu02}O. R. Musin, The kissing number in four dimensions, \emph{Ann. of Math.} \textbf{168} (2008), 1-32.
\bibitem{Od01}A. M. Odlyzko, N. J. A. Sloane, New bounds on the number of unit spheres that can touch a unit sphere in $n$ dimensions, \emph{J. Combin. Theory Ser. A} \textbf{26} (1979), 210-214.
\bibitem{Sc02}A. Sch\"{u}rmann, F. Vallentin, Computational Approaches to Lattice Packing and Covering Problems, \emph{Discrete Comput Geom.} \textbf{35} (2006), 73-116.
\bibitem{Sc01}K. Sch\"{u}tte, B. L. van der Waerden, Das Problem der dreizehn Kugel, \emph{Math. Ann.} \textbf{125} (1953), 325-334.
\bibitem{Sw01}H. P. F. Swinnerton-Dyer, Extremal lattices of convex bodies, \emph{Math. Proc. Cambridge Philos. Soc.} \textbf{49} (1953), 161-162.
\bibitem{Wa01}G. L. Watson, The number of minimum points of a positive quadritic form, \emph{Dissertationes Math.} \textbf{84} (1971), 1-42.
\bibitem{Zo02}C. Zong, The kissing numbers of tetrahedra, \emph{Discrete Comput. Geom.} \textbf{15} (1996), 239-252.
\bibitem{Zo03}C. Zong, \emph{Strange Phenomena in Convex and Discrete Geometry}, Springer-Verlag, New York, 1996.
\bibitem{Zo06}C. Zong, \emph{Sphere Packings}, Springer-Verlag, New York, 1999.
\bibitem{Zo04}C. Zong, From deep holes to free planes, \emph{Bull. Amer. Math. Soc.} \textbf{39} (2002), 533-555.
\bibitem{Zo07}C. Zong, Simultaneous packing and covering in three-dimensional Euclidean space, \emph{J. London Math. Soc.} \textbf{67} (2003), 29-40.
\end{thebibliography}
\end{document}